\documentclass[
journal=jacsat, 
manuscript=article]{achemso}

\usepackage[version=3]{mhchem} 

\SectionNumbersOn
\usepackage{graphicx}
\usepackage[english]{babel}
\usepackage{wrapfig}
\usepackage[T1]{fontenc}
\usepackage{mathptmx}
\usepackage[english]{babel}
\usepackage{amsmath}
\usepackage{mathtools}
\usepackage{amsfonts}
\usepackage{amssymb}
\usepackage{makeidx}
\usepackage{graphicx}
\usepackage{scalerel}
\usepackage{hyperref}
\usepackage{placeins}
\usepackage{xfrac}
\usepackage{color}
\usepackage{float}
\usepackage{caption}
\usepackage{subcaption}
\usepackage{chngcntr}
\usepackage{setspace}
\usepackage{array}
\usepackage[labelsep=space]{caption}
 \setcitestyle{aysep={}} 
\author{Hammed O. Fatoyinbo}
\email{h.fatoyinbo@massey.ac.nz}
\author{Richard G. Brown}
\author{David J. W. Simpson}
\author{Bruce van Brunt}
\affiliation[Massey University]
{School of Fundamental Sciences, Massey University, New Zealand}

\title{Numerical Bifurcation Analysis of Pacemaker Dynamics in a Model of Smooth Muscle Cells}

\begin{document}
\begin{abstract}\label{abstract}
Evidence from experimental studies shows that oscillations due to electro-mechanical coupling  can be generated spontaneously in smooth muscle cells. Such cellular dynamics are known as \textit{pacemaker dynamics}. In this article we address pacemaker dynamics associated with the interaction of $\text{Ca}^{2+}$ and $\text{K}^+$ fluxes in the cell membrane of a smooth muscle cell. First we reduce a pacemaker model to a two-dimensional system equivalent to the reduced Morris-Lecar model and then perform a detailed numerical bifurcation analysis of the reduced model. Existing bifurcation analyses of the Morris-Lecar model concentrate on external applied current whereas we focus on parameters that model the response of the cell to changes in transmural pressure. We reveal a transition between Type I and Type II excitabilities with no external current required. We also compute a two-parameter bifurcation diagram and show how the transition is explained by the bifurcation structure.
\end{abstract}

 \newpage

\section{Introduction}

\label{sec:introduction}
Electro-mechanical coupling (EMC) refers to the contraction of a smooth muscle cell (SMC) due to its excitation in response to an external mechanical stimulation, such as a change in transmural pressure, that is, the pressure gradient across the vessel wall \citep{Ran2019TransmuralSurvivin}.  In some SMCs, EMC activity can be spontaneous owing to interactions between ion fluxes through voltage-gated ion channels. Based on experimental observations, c.f. \citep{ Casteels1977Excitation-contractionArtery.,Harder1984Pressure-DependentArtery,Koenigsberger2005RoleVasomotion}, the ion channels coordinating the EMC activity in SMCs of feline cerebral arteries are the voltage-gated  $\text{Ca}^{2+}$ channel, voltage-gated $\text{K}^{+}$ ion channel and the leak ion channel. The spontaneous depolarisation of the cell membrane leads to the opening and closing of ion channels resulting in a fluctuation of ionic currents that can induce EMC activity \citep{Sui2003AMuscle,  Brading2006SpontaneousFunction, Mahapatra2018AMuscle}. This pacemaker EMC activity varies across species of SMCs \citep{Savineau2000CytosolicCells} and understanding the impact  of these dynamics on the type of excitability may suggest therapeutic strategies for treating diseases related to SMCs. 

Under normal physiological conditions, the cell membranes of SMCs do not oscillate in the absence of external sources, however several exceptions have been observed. \citet{Mclean1977} studied the spontaneous contraction of cultured vascular SMCs in chick embryos.  \cite{Lusamvuku1979CorrelatedDrugs} observed spontaneous electrical activity in rabbit cerebral arteries exposed to high pressure.  \cite{Harder1984Pressure-DependentArtery} examined  cellular  mechanisms  of  the  myogenic response, the pressure-induced contraction of blood vessels to regulate blood flow, in feline middle  cerebral arteries by recording  intracellular electrical activity of arterial muscle cells upon elevation  of  transmural pressure. It was observed that the blood vessels contract and spontaneous firing occurs as the arterial blood pressure is increased. \citet{Llinas1988TheFunction} experimentally explored auto-rhythmic electrical properties in the mammalian central nervous system. \citet{Meister1991} and \citet{Gu2013ExperimentalPacemaker} reported experimental observations of spontaneous oscillations induced by modulating either extracellular calcium or potassium concentrations in neural cells.

Research into EMC activity  has shown that abnormal contraction is often associated with tissue diseases. For example abnormal vasomotion in arteries can damage blood vessels causing hypertension over time \citep{Humphrey2003AStudy} and spontaneous contraction of the urinary bladder causes urine leak \citep{Brading2006SpontaneousFunction}.

The dynamics of electrical activity in cell membranes are nonlinear, and often well-modelled by a nonlinear system of ODEs \citep{Izhikevich2007DynamicalBursting, Ma2015ANeurons}. Many such models have been developed to describe the behaviour of excitable cells in the cell membrane. The pioneering work of Hodgkin and Huxley describes the conduction of electrical impulses along a squid giant axon \citep{Hodgkin1952}. Other well known models include the FitzHugh-Nagumo model \citep{Fitzhugh, Nagumo1962}, the Morris-Lecar model \citep{Morris}, the Hindmarsh-Rose model \citep{Indmarsh11984}, and the Izhikevich model \citep{Izhikevich2007DynamicalBursting}.  

As revealed in experiments, the electrical activity of a single excitable cell has a variety of possible dynamical behaviours, such as a rest or quiescent state,
simple oscillatory motion, and complex oscillatory motion. A model of a cell can transition from one state to another as parameters are varied \citep{Gu2013BiologicalModel,Gu2013ExperimentalPacemaker}. These changes can be understood by identifying critical parameter values (bifurcations) at which the dynamical behaviour  changes qualitatively \citep{Strogatz1994NonlinearEngineering, KuznetsovY.A.1995ElementsTheory, Meiss2007}. For excitable cells, arguably the  most important transition is from rest to an oscillatory state (or vice versa). Bifurcations associated with this and other transitions have been identified in many studies \citep{Govaerts2005TheApproach, Tsumoto2006BifurcationsModel, Prescott2008BiophysicalInitiation, Storace2008TheApproximations, Barnett2014AModel, Liu2014BifurcationModel, Zhao2017TransitionsAutapse,  Mondal2018DynamicsSolutions, Mondal2019BifurcationModel}.  

Models for excitable cells can be classified into two types depending on the nature of action potential generation. \cite{Rinzel1999AnalysisNetwork} used the type of bifurcation at the onset of firing to classify excitable cells into Type I and Type II. In Type I excitability, the cell transitions from rest to an oscillatory state through a saddle-node on an invariant circle (SNIC) bifurcation. As parameters are varied to move away from the bifurcation, the frequency of the oscillations increases from zero. In contrast, for Type II excitability the transition from rest to an oscillatory state is through a Hopf bifurcation.  In this case the oscillations emerge with non-zero frequency. Rinzel and Ermentrout \citep{Rinzel1999AnalysisNetwork} also concluded that  their classification is consistent with the original classification of \cite{Hodgkin1948TheAxon} for the squid giant axon, see \citep{Ermentout1996TypeSychrony, Rinzel1999AnalysisNetwork, Crook1998SpikeOscillators, Vreeswijk2001PatternsAdaptation}

 The Morris-Lecar model can exhibit both Type I and Type II excitability depending on the parameter regime. \cite{Rinzel1999AnalysisNetwork} studied Type I and Type II  excitability in the reduced Morris-Lecar model by adjusting the applied current. \cite{Tsumoto2006BifurcationsModel} and \cite{Zhao2017TransitionsAutapse} subsequently identified codimension-two bifurcations associated with a change between the two types of excitability. See also \cite{Duan2008Two-parameterModel} for a similar two-parameter bifurcation analysis of the Chay neuronal model. 

Recently there have been several studies of pacemaker dynamics in excitable cells, both theoretical \citep{Duan2006Codimension-twoModel,Duan2008Two-parameterModel} and computational \citep{Gonzalez-Miranda2012NonlinearCell}. The importance of the leak channel in the pacemaker dynamics of the full Morris-Lecar model has been studied by  \citet{Gonzalez-Miranda2014}. Also \citet{Meier} confirmed the existence of spontaneous action potentials in the two variable Morris-Lecar model.  Despite many studies of pacemaker activity in SMCs having being conducted, there does not appear to have been any discussion about the types of excitabilility that can be exhibited.

The purpose of this paper is to explain the occurrence of Type I and II excitability in pacemaker dynamics. We begin in Sect.~\ref{sec:model} with the three-dimensional ODE model of \cite{Gonzalez-Fernandez1994} for pacemaker dynamics in feline cerebral arteries. In Sect.~\ref{sec:model_reduction} we apply a small simplification to the model which reduces the dynamics to the two-variable Morris-Lecar model with no applied current and nondimensionalise the model in Sect.~\ref{sec:nondimensionalised}.

Then in Sect.~\ref{sec:bifurcation} we perform a detailed bifurcation analysis of the nondimensionalised model. As the primary bifurcation parameter we use the voltage associated with the opening of the $\text{K}^+$ channels because experiments have revealed that action potentials can be triggered by an increase in transmural pressure  \citep{Harder1984Pressure-DependentArtery,Harder}. We find both types of excitability and identify codimension-two bifurcations that represent endpoints for the two types of excitability. We stress that while the bifurcations we find have been described already in the Morris-Lecar model \citep{Tsumoto2006BifurcationsModel, Zhao2017TransitionsAutapse}, we believe that this is the first work to describe this structure in pacemaker dynamics of SMCs. Moreover this work is a necessary first step towards understanding spatiotemporal behaviour in networks of SMCs connected electrically by gap junctions. Finally conclusions are presented in Sect.~\ref{sec:conclusion}

\newpage
\section{Model Formulation}\label{sec:model_equation}

  \subsection{Muscle Cell Model}\label{sec:model}
 \cite{Gonzalez-Fernandez1994}  consider a muscle cell model with external current set to zero to study pacemaker dynamics. The model consists of the three ODEs
 
\begin{align}
\label{eq:V}
\text{C}\frac{dv}{dt}&=-g_{L}(v-v_{L})-g_{K}n(v-v_{K})-g_\text{Ca}m_{\infty}(v)(v-v_\text{Ca}),\\
\label{eq:N}
\frac{dn}{dt}&=\lambda_{n}(v)\big(n_{\infty}(v,\text{Ca}_{i})-n\big),\\
\label{eq:Ca}
\frac{d\text{Ca}_{i}}{dt}&=\big(-\alpha g_\text{Ca}m_{\infty}(v)(v-v_\text{Ca})-k_\text{Ca}\text{Ca}_{i}\big)\rho(\text{Ca}_i),
\end{align}

\noindent where $v$ is the membrane potential, $n$ is the fraction of open potassium channels, and $\mbox{Ca}_i$ is the cytosolic concentration of calcium. The system parameters $g_{L}$, $g_{K}$, and $g_\text{Ca}$ are the maximum conductances for the leak, potassium, and calcium currents, respectively, while $v_{L}$, $v_{K}$ and  $v_\text{Ca}$ are the corresponding Nernst reversal potentials. Also $\text{C}$ is the cell capacitance, $k_\text{Ca}$ is the rate constant for cytosolic calcium concentration, and $\rho$ models the calcium buffering.  The auxiliary functions in the model are:

\begin{align}
\label{eq:MINF}
m_\infty(v)&=0.5\left(1+\tanh\left(\frac{v-v_{1}}{v_{2}}\right)\right),\\
\label{eq:nif}
n_{\infty}(v,\text{Ca}_{i})&=0.5\left(1+\tanh\left(\frac{v-v_{3}(\text{Ca}_i)}{v_{4}}\right)\right),\\
\label{eq:v3}
v_{3}(\text{Ca}_i)&=-\frac{v_{5}}{2}\tanh\left(\frac{\text{Ca}_{i}-\text{Ca}_{3}}{Ca_{4}}\right)+v_{6},\\
\label{eq:lam}
\lambda_{n}(v)&=\phi_{n}\cosh\left(\frac{v-v_{3}(\text{Ca}_i)}{2v_{4}}\right),\\
\label{eq:rho}
\rho(\text{Ca}_{i})&=\frac{(K_{d}+\text{Ca}_{i})^{2}}{(K_{d}+\text{Ca}_{i})^{2}+K_{d}B_{T}},
\end{align}

\noindent where $n_{\infty}$ [$m_\infty$] is the fraction of open potassium [calcium] channels at steady state, $\phi_{n}$ is the rate constant for the kinetics of the potassium channel, $K_{d}$ is the ratio of backward and forward binding rates for calcium and buffer reaction \citep{Sala1990}, and $B_{T}$ is the total concentration of the buffers. For further details see \cite{Gonzalez-Fernandez1994}. The  parameter values of  \cite{Gonzalez-Fernandez1994} are listed in Table \ref{Tab:1}.

\begin{table}[h!]
 \centering
 \caption{Model parameter values are taken from \cite{Gonzalez-Fernandez1994}} 
\begin{tabular}{lll}
\hline\noalign{\smallskip}
Parameter &Value& Unit\\ 
\hline
$v_{1}$& $-22.5$ & mV \\ 

$v_{2}$  & $25.0$  & mV\\ 

$v_{4}$ & $14.5$  & mV \\ 

$v_{5}$ & $8.0$ & mV \\ 

$v_{6}$& $-15.0$  & mV \\ 

$\text{Ca}_{3}$& $400.0$ & nM \\ 

$\text{Ca}_{4}$& $150.0$  & nM\\ 

$\phi_{n}$  & $2.664$ & $s^{-1}$ \\ 

$v_{L}$ & $-70.0$  & mV  \\ 

$v_{K}$   & $-90.0$ & mV \\ 

$v_{\text{Ca}}$   & $80.0$& mV  \\ 

$\text{C}$ & $1.9635\times 10^{-14}$ & Cm$V^{-1}$ \\ 

$g_{L}$ & $7.854\times 10^{-14}$  & C$s^{-1}$m$V^{-1}$\\ 

$g_{K}$  & $3.1416\times 10^{-13} $ & C$s^{-1}$m$V^{-1}$\\ 

$g_{\text{Ca}}$  & $1.57\times 10^{-13}$ & C$s^{-1}$m$V^{-1}$ \\ 

$K_{d}$& $1.0\times 10^{3}$ & nM  \\ 

$B_{T}$ & $1.0\times 10^{5}$  & nM \\ 

$\alpha$& $7.9976\times 10^{15}$  & nM$\text{C}^{-1}$ \\ 

$k_{\text{Ca}}$  & $1.3567537\times 10^{2}$ & $s^{-1}$ \\ 
\noalign{\smallskip}\hline 
\end{tabular} 
\label{Tab:1}
\end{table}
 
\subsection{Model Reduction} \label{sec:model_reduction}
\begin{figure}
\centering
\begin{subfigure}[b]{.5\textwidth}
    \centering
    \caption{}
    \includegraphics[width = \textwidth]{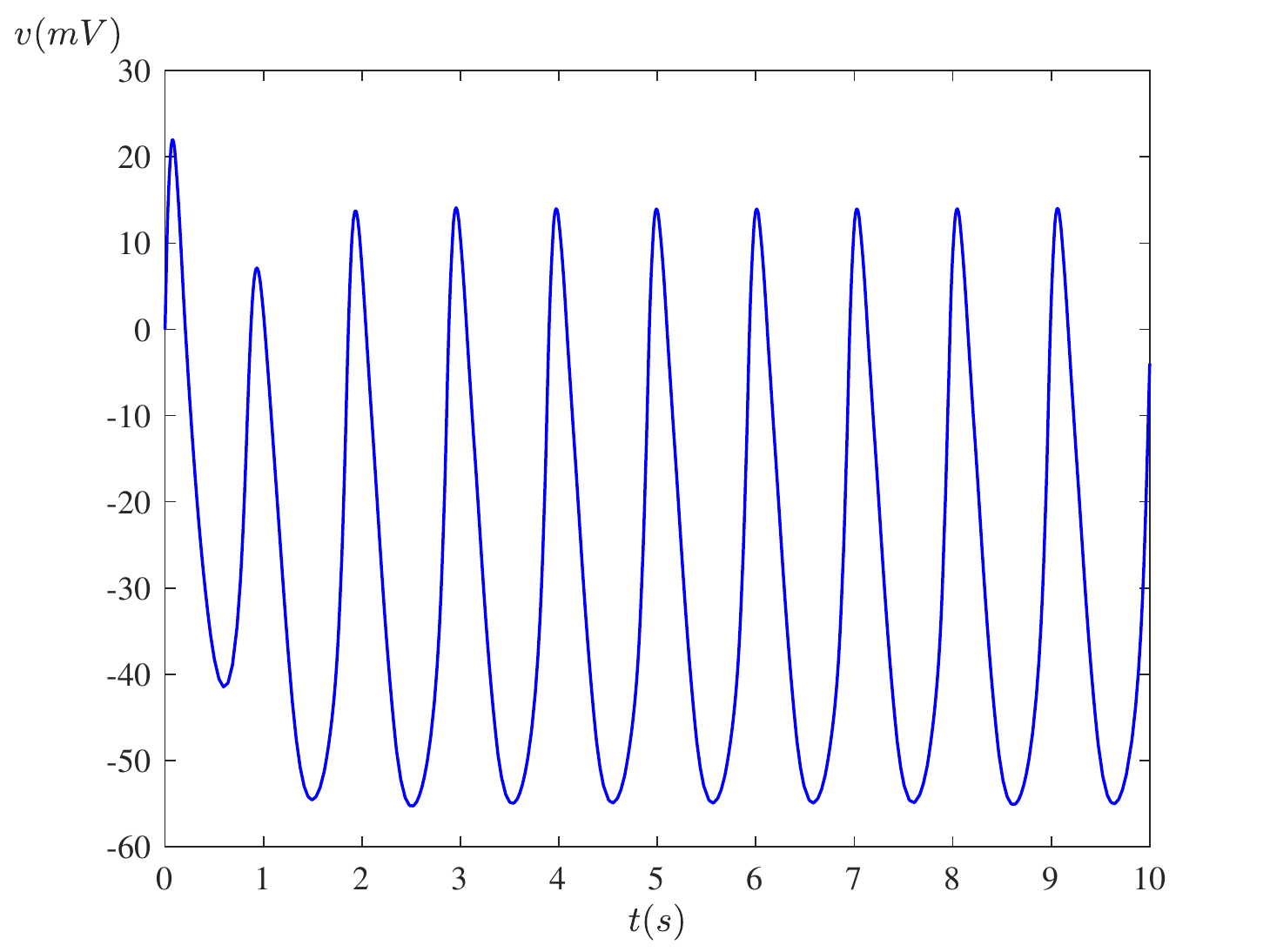}
  \end{subfigure}\\%
  \begin{subfigure}[b]{.5\textwidth}
    \centering
    \caption{}
    \includegraphics[width = \textwidth]{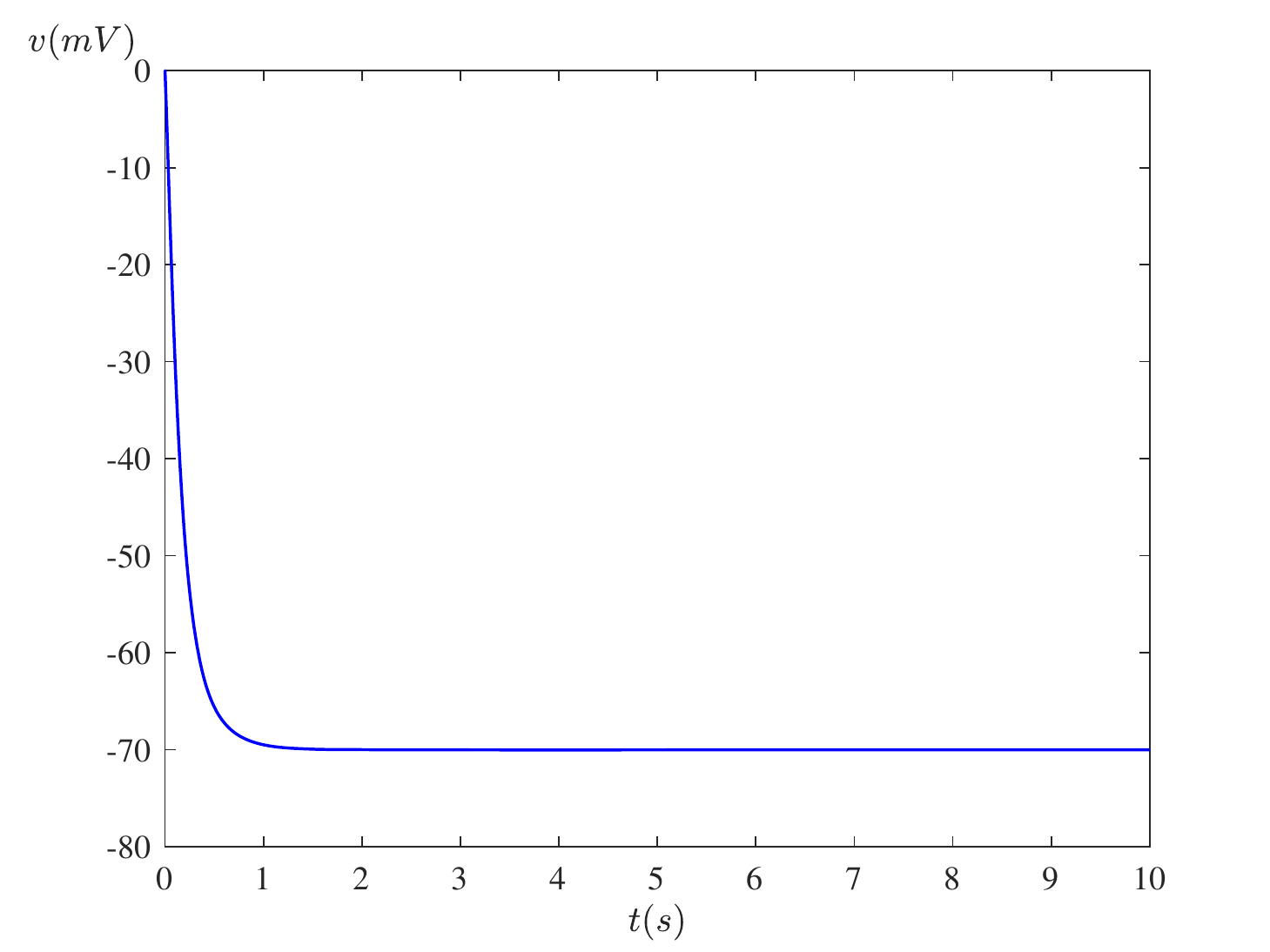}
  \end{subfigure}%
  \begin{subfigure}[b]{.5\textwidth}
    \centering
    \caption{}
    \includegraphics[width = \textwidth]{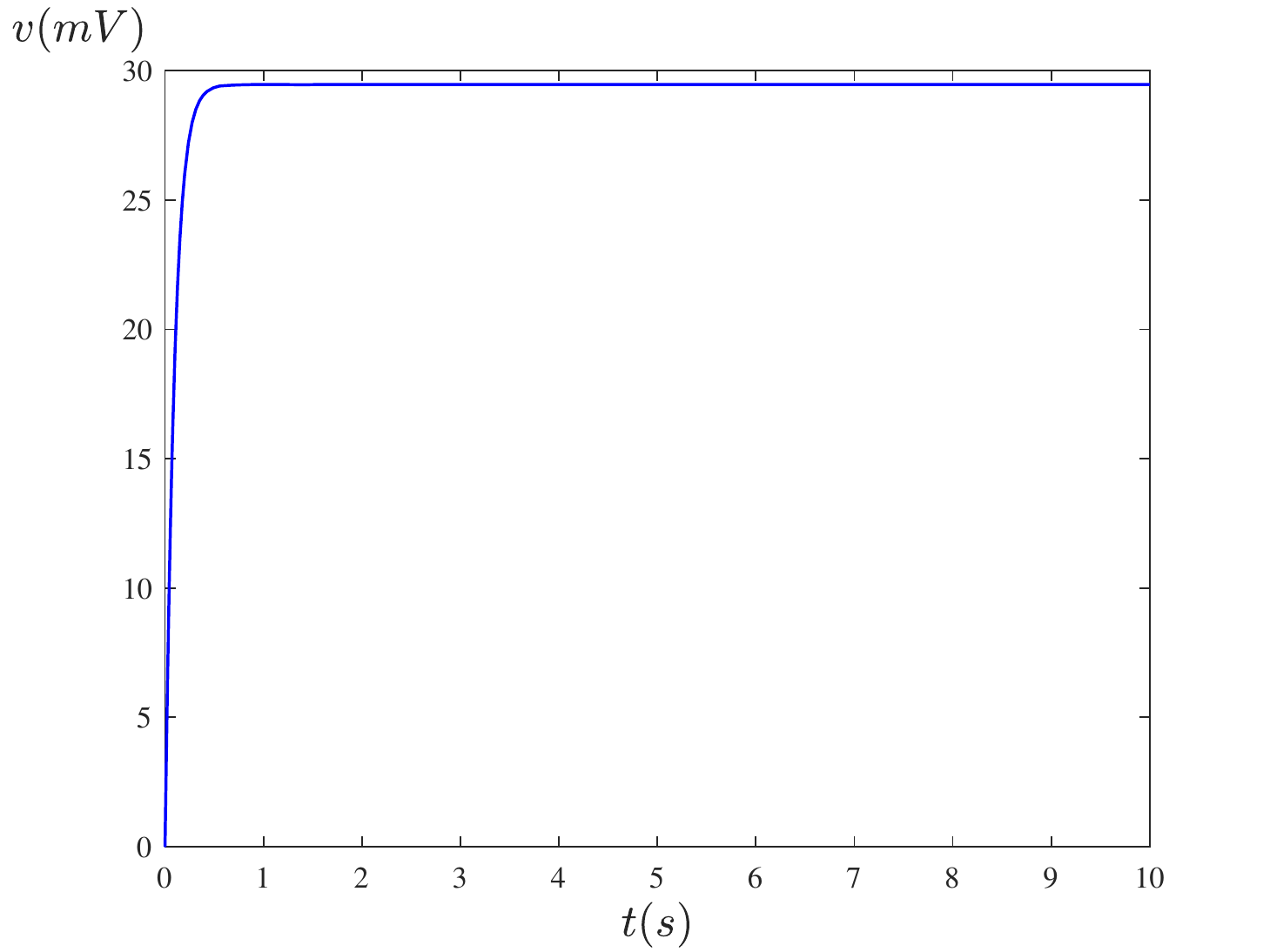}
  \end{subfigure}%
 \caption{ Time series of the membrane potential $v$ when the three conductances are blocked: (a) the leak channel is blocked ($g_\text{L}$); (b) the $\text{Ca}^{2+}$ channel is blocked ($g_\text{Ca}$); (c) the $\text{K}^{+}$ channel is blocked ($g_\text{K}$)}
\label{fig:1}
\end{figure}

 To analyse the model we first check the effects of each ionic current on pacemaker activity. To do this we block the conductances for the leak, $\text{Ca}^{2+}$, and $\text{K}^{+}$ currents in turn. Over a range of parameter values we found that pacemaker activity persists if the leak current conductance $g_L$ is blocked, but is absent if the conductances  $g_{\text{Ca}}$ and $g_K$ for the $\text{Ca}^{2+}$ and $\text{K}^{+}$  currents are blocked (Fig.~\ref{fig:1} shows an example). This tells us that the $\text{Ca}^{2+}$ and $\text{K}^{+}$  currents are required for pacemaker activity in the model.

\begin{figure}
\centering
\begin{subfigure}{.6\textwidth}
  \centering
  \includegraphics[width=\textwidth]{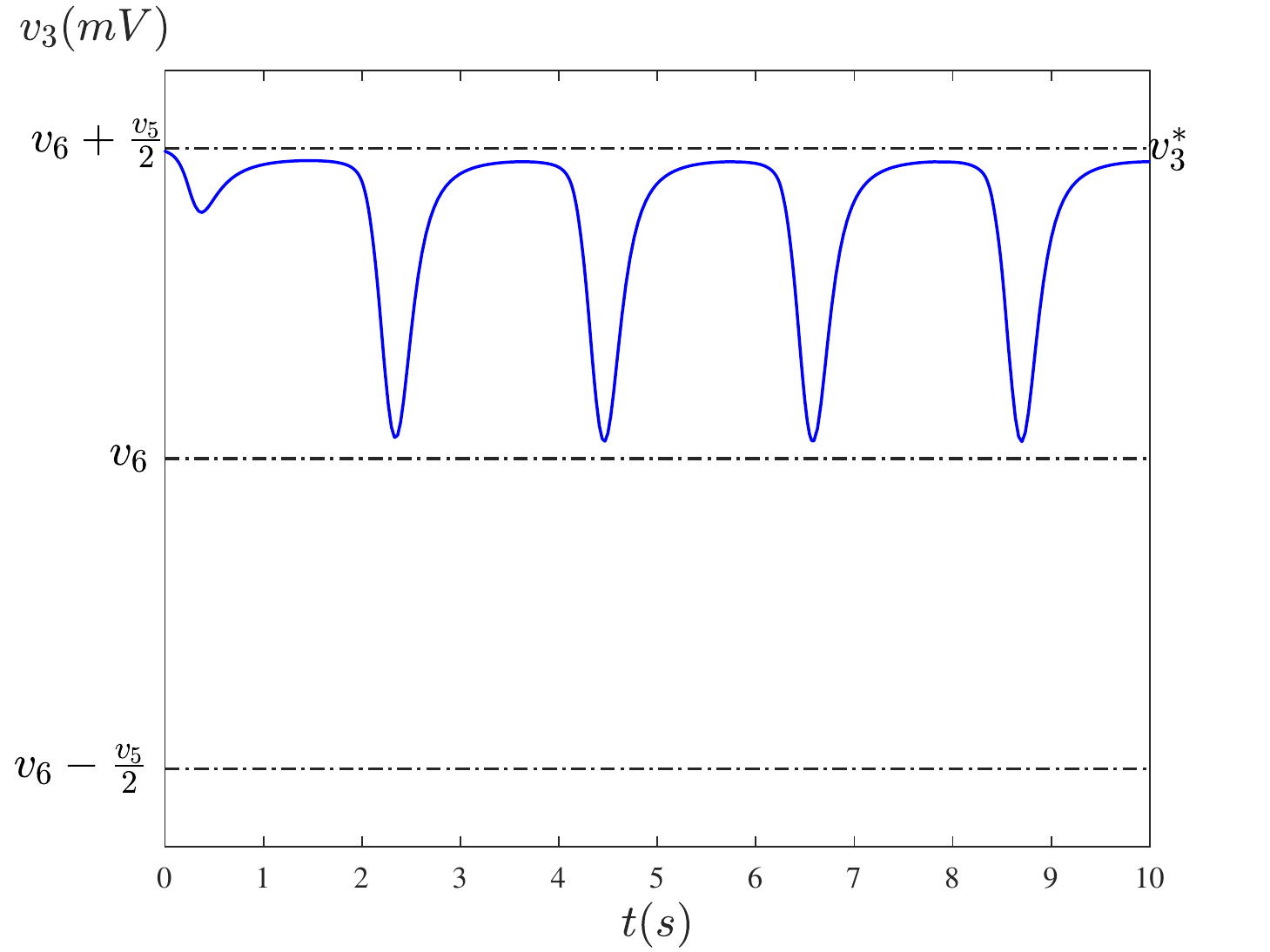}
\end{subfigure}
\caption{ A plot of $v_3(mV)$ against time for solutions to \eqref{eq:V}--\eqref{eq:Ca} with the parameters of Table \ref{Tab:1}}
\label{fig:2}
\end{figure}

We now reduce system \eqref{eq:V}--\eqref{eq:Ca} to two equations. Our reduction is based on the behaviour of the time-dependent quantity $v_3$. Equation \eqref{eq:v3} shows that the value of $v_3$ has the upper and lower bounds $v_6 + \frac{v_5}{2}$ and $v_6 - \frac{v_5}{2}$, respectively. Using the parameter values of Table \ref{Tab:2} and a numerical solution to system \eqref{eq:V}--\eqref{eq:Ca}, we see from Fig.~\ref{fig:2} that the value of $v_3$ spends a high proportion of time close to its upper bound (after transient dynamics have decayed). This motivates a reduction by fixing $v_3$ to the value of its upper bound. We thus replace \eqref{eq:v3} with $v_3 = v_3^*$, where $v_3^* = v_6 + \frac{v_5}{2}$. See already Fig.~\ref{fig:4} which shows that the bifurcation structure of the resulting reduced model is similar to that of the full model. The equilibria undergo the same sequences of bifurcations in the same order, which indicates that the reduction does not significantly alter the qualitative dynamics. The assumption of constant $v_3$ reduces the number of equations to two because now $v$ and $n$ are decoupled from $\text{Ca}_i$. The reduced system is 

\begin{align}
\label{eq:reduced_2}
\text{C}\frac{dv}{dt}&=-g_{L}(v-v_{L})-g_{K}n(v-v_{K})-g_\text{ca}m_{\infty}(v)(v-v_\text{ca}),\\
\label{eq:reduced_3}
\frac{dn}{dt}&=\lambda_{n}(v)\left(n_{\infty}(v)-n\right),
\end{align}
where
\begin{align}
\label{eq:ninf}
    n_{\infty}(v)&=0.5\left(1+\tanh\left(\frac{v-v_{3}^*}{v_{4}}\right)\right),\\
\label{eq:lam_1}
\lambda_{n}(v)&=\phi_{n}\cosh\left(\frac{v-v_{3}^*}{2v_{4}}\right),
\end{align}
and $m_\infty(v)$ is unchanged from \eqref{eq:MINF}. Note that this is the Morris-Lecar model without external current.

\newpage
\subsection{Nondimensionalised model}\label{sec:nondimensionalised}
We nondimensionalise \eqref{eq:reduced_2}--\eqref{eq:reduced_3} by introducing dimensionless variables $V$ and $\tau$. Let
\begin{equation}
    v=VQ_v, \hspace{0.5cm} t=\tau Q_t,
    \label{eq:transformation}
\end{equation}
for some characteristic voltage $Q_v$ and  time $Q_t$. To choose values for $Q_v$ and $Q_t$ we first observe that the range of the action potential is $v_K\leq v \leq v_\text{Ca}$ (see Table \ref{Tab:1} and Fig.~\ref{fig:2.3a}). Hence the maximum variation of the action potential is less than $v_\text{Ca}-v_{K}=170$mV. This value is roughly the same order of magnitude as $v_{\text{Ca}}$ therefore we choose the characteristic voltage $Q_v$ to  be $v_\text{Ca}$. Simple choices for the characteristic time include $Q_{t} = \frac{C}{g_K} = 0.0625$ and $Q_{t}=\frac{1}{\phi_n}=0.3754$. We choose $Q_{t}=\frac{C}{g_K}$ for the characteristic time because it is faster than $\frac{1}{\phi_n}$. Substituting $Q_v=v_\text{Ca}$ and $Q_t=\frac{C}{g_K}$ into \eqref{eq:reduced_2}--\eqref{eq:reduced_3} produces the dimensionless version of the model:
\begin{align}
\label{eq:dimless1}
\frac{dV}{dT}&=-\bar{g}_{L}(V-\bar{v}_{L})-\bar{g}_{K}N(V-\bar{v}_{K})-\bar{g}_{\text{Ca}}M_{\infty}(V)(V-1), \\
\frac{dN}{dT}&=\psi\lambda(V)(N_{\infty}(V)-N),
\label{eq:dimless2}
\end{align}
where
\begin{align}
\label{MINF_2}
M_{\infty}(V)&=0.5\left(1+\tanh\left(\frac{V-\bar{v}_{1}}{\bar{v}_{2}}\right)\right),\\
N_{\infty}(V)&=0.5\left(1+\tanh\left(\frac{V-\bar{v}_{3}}{\bar{v}_{4}}\right)\right),\\
\lambda(V)&=\cosh\left(\frac{V-\bar{v}_{3}}{2\bar{v}_{4}}\right),
\label{eq:lambda}
\end{align}
and $$\bar{g}_i=\frac{g_i}{g_K}, \hspace{0.2cm} \bar{v}_i=\frac{v_i}{v_\text{Ca}}, \hspace{0.2cm} \psi=\frac{C\phi_n}{g_K}, \hspace{1cm} i=L,K,\text{Ca}, 1,2,3,4.$$ 
The parameter values for this model are given in Table \ref{Tab:2}. 

\begin{table}[H]
\centering
\caption{Parameter values for the  nondimensionalised model \eqref{eq:dimless1}--\eqref{eq:dimless2}}
\label{Tab:2}
\begin{tabular}{ll}
\hline\noalign{\smallskip}
Parameter& Value\\
\hline
$\bar{v}_{1}$ &  $-0.2813$ \\ 

$\bar{v}_{2}$ &  $0.3125$ \\ 

$\bar{v}_{3}$& $-0.1380$ \\ 

$\bar{v}_{4}$ & $-0.1812$ \\

$\psi$ & $0.1665$\\

$\bar{v}_{L}$ &  $-0.875$  \\ 

$\bar{v}_{K}$  &  $-1.125$  \\ 

$\bar{g}_{L}$ & $0.25$ \\ 

$\bar{g}_{K}$ &  $1.0 $\\ 

$\bar{g}_{\text{Ca}}$  & $0.4997$ \\ 
\noalign{\smallskip}\hline
\end{tabular} 
\end{table}

\subsection{Excitable dynamics of the full, reduced, and nondimensionalised models.}\label{sec:dynamics}

\begin{figure}[!ht]
\centering
\begin{subfigure}{.5\textwidth}
 \centering
 \caption{}
  \includegraphics[width=\textwidth]{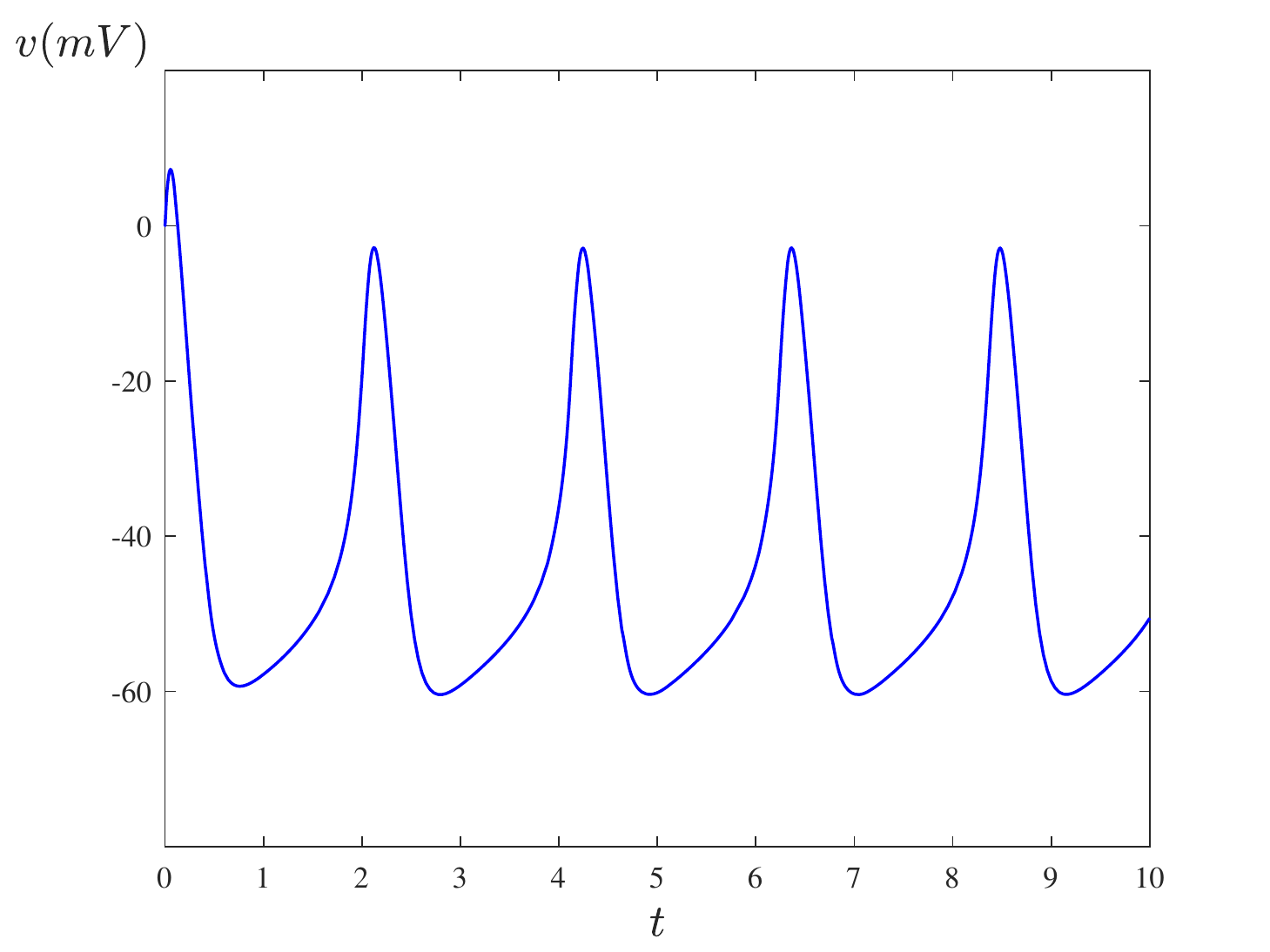}
  \label{fig:2.3a}
\end{subfigure}\\%
\begin{subfigure}{.5\textwidth}
  \centering
  \caption{}
  \includegraphics[width=\textwidth]{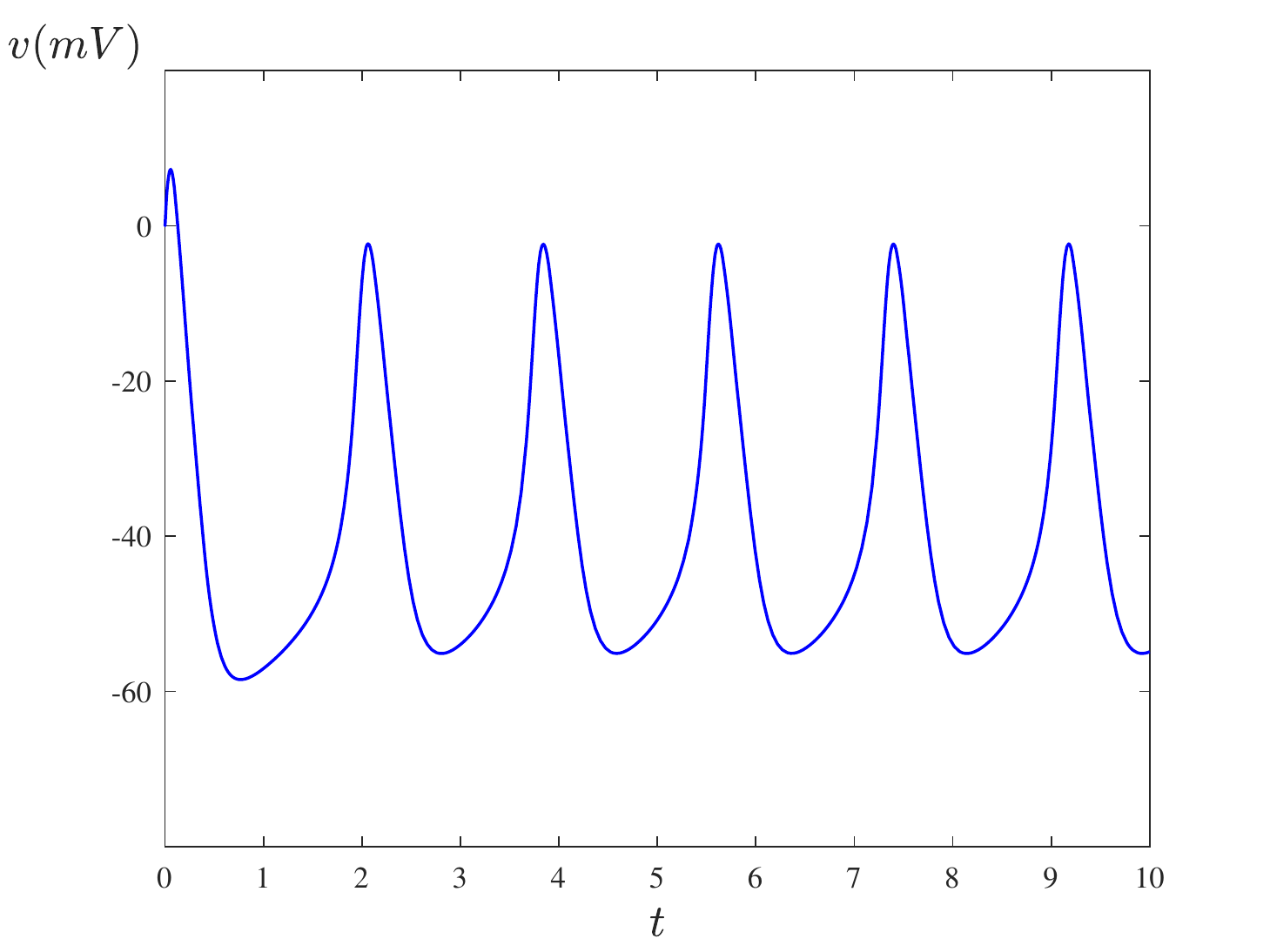}
\label{fig:2.3b}
\end{subfigure}%
\begin{subfigure}{.5\textwidth}
     \centering
     \caption{}
     \includegraphics[width=\textwidth]{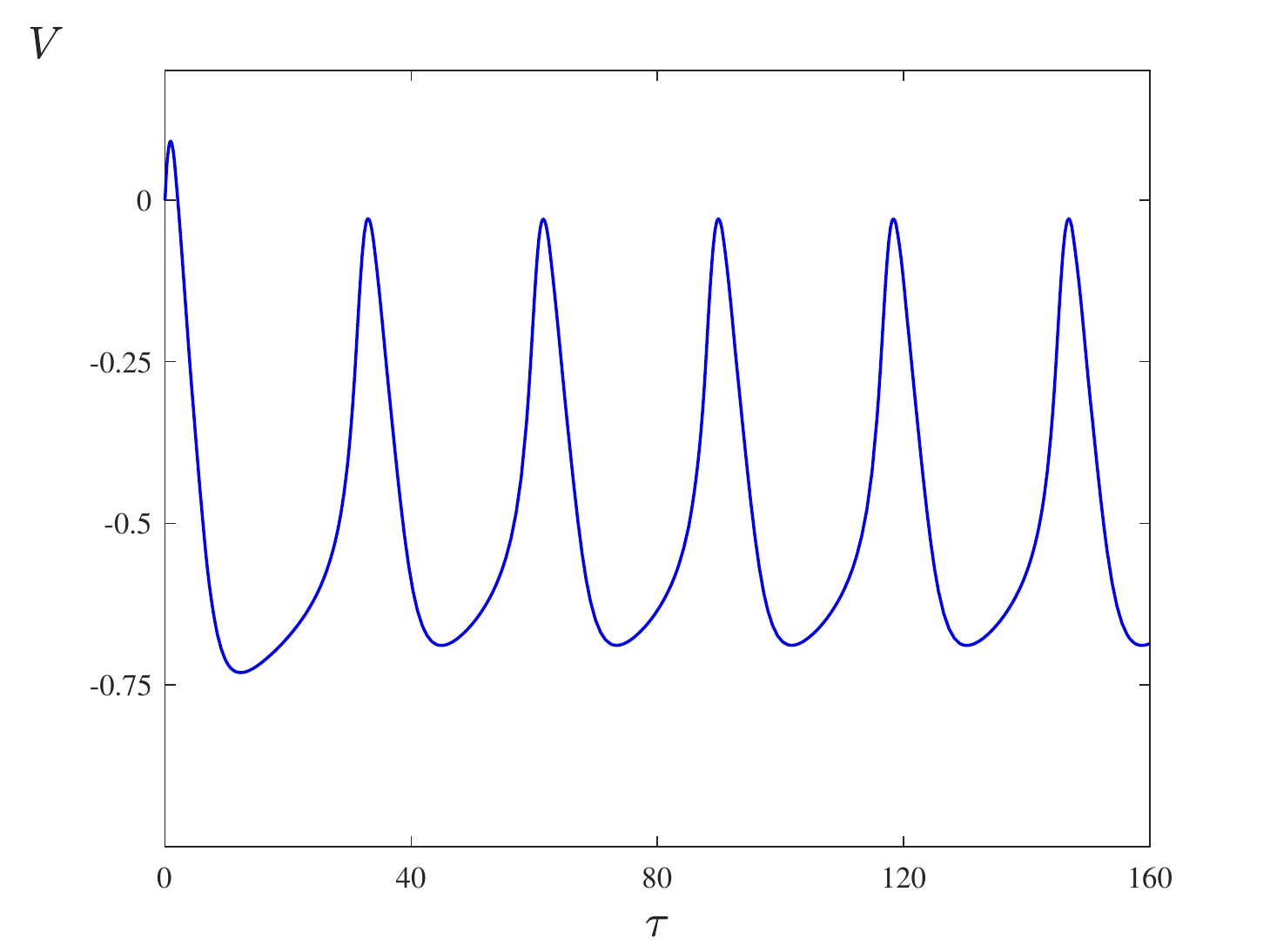}
     \label{fig:2.3c}
      \end{subfigure}
      \caption{A time series of the membrane potential for (a) the full model with the parameter values in Table \ref{Tab:1}  and initial condition $(v,n,\text{Ca}_i)=(0,0,0)$ (b) the reduced model , and (c) the nondimensionalised model with the parameter values in Table \ref{Tab:2} and initial condition $(V,N)=(0,0)$ }
\label{fig:3}
\end{figure}

The full  model \eqref{eq:V}--\eqref{eq:Ca}, the reduced model \eqref{eq:reduced_2}--\eqref{eq:reduced_3}, and the nondimensionalised model \eqref{eq:dimless1}--\eqref{eq:dimless2}  were integrated numerically using the standard fourth-order Runge-Kutta method using a step size of 0.05 in the numerical software XPPAUT \citep{Ermentrout2002SimulatingStudents}. Since our interest is primarily the membrane potential, we focus mostly on its dynamics. The time evolution of the membrane potential for the three models with the parameter values in Tables \ref{Tab:1} and \ref{Tab:2} reveal that they are in an oscillatory state (see Fig. \ref{fig:3}). These self-sustained oscillations are consistent with the work of \cite{Gonzalez-Miranda2014} on pacemaker dynamics for the full Morris-Lecar model when the external current and the leak conductance are set to zero. 

\begin{figure}
\centering
\begin{subfigure}{.5\textwidth}
 \centering
 \caption{}
  \includegraphics[width=\textwidth]{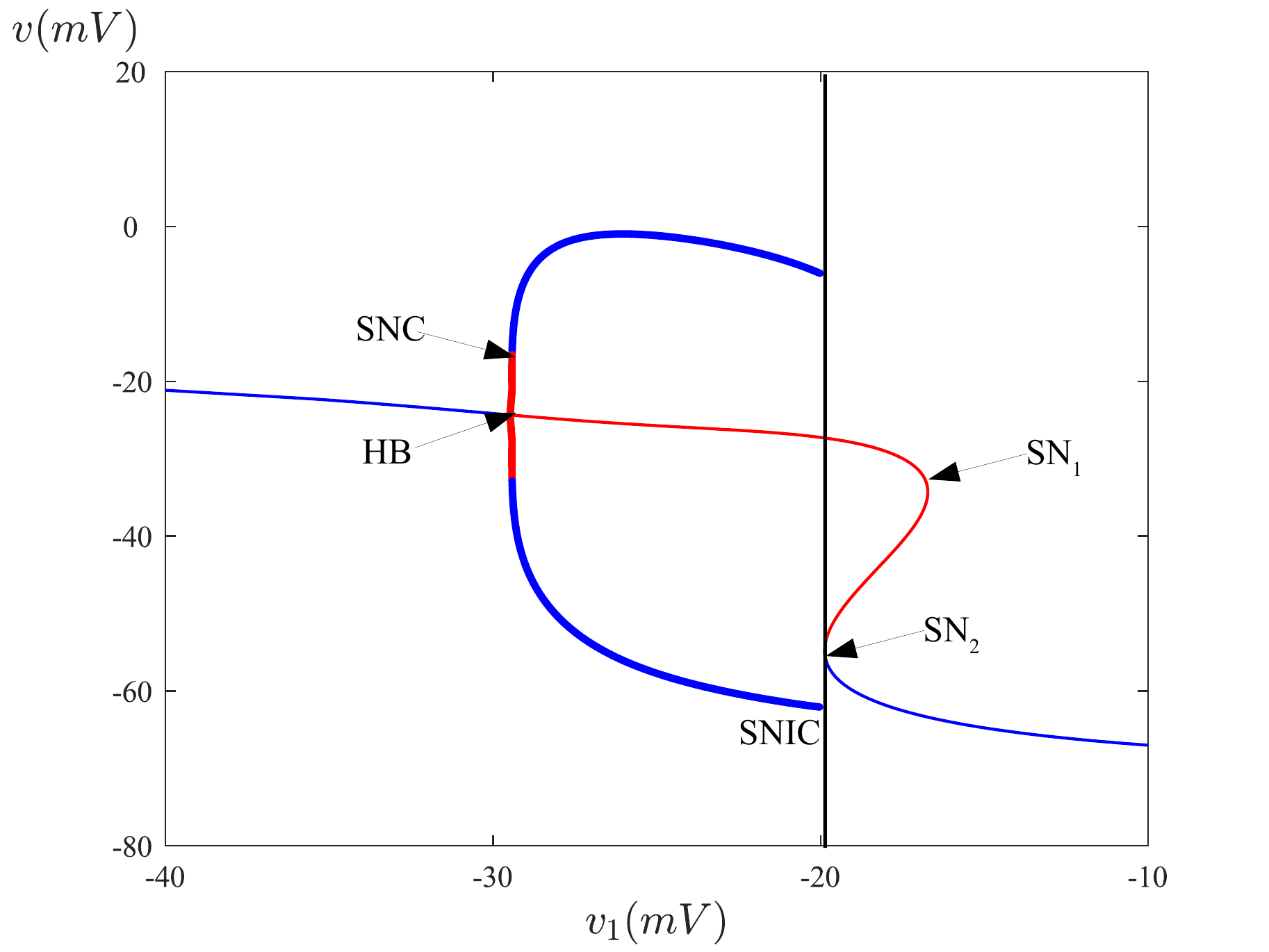}
  \label{fig:2.4a}
\end{subfigure}%
\begin{subfigure}{.5\textwidth}
  \centering
  \caption{}
  \includegraphics[width=\textwidth]{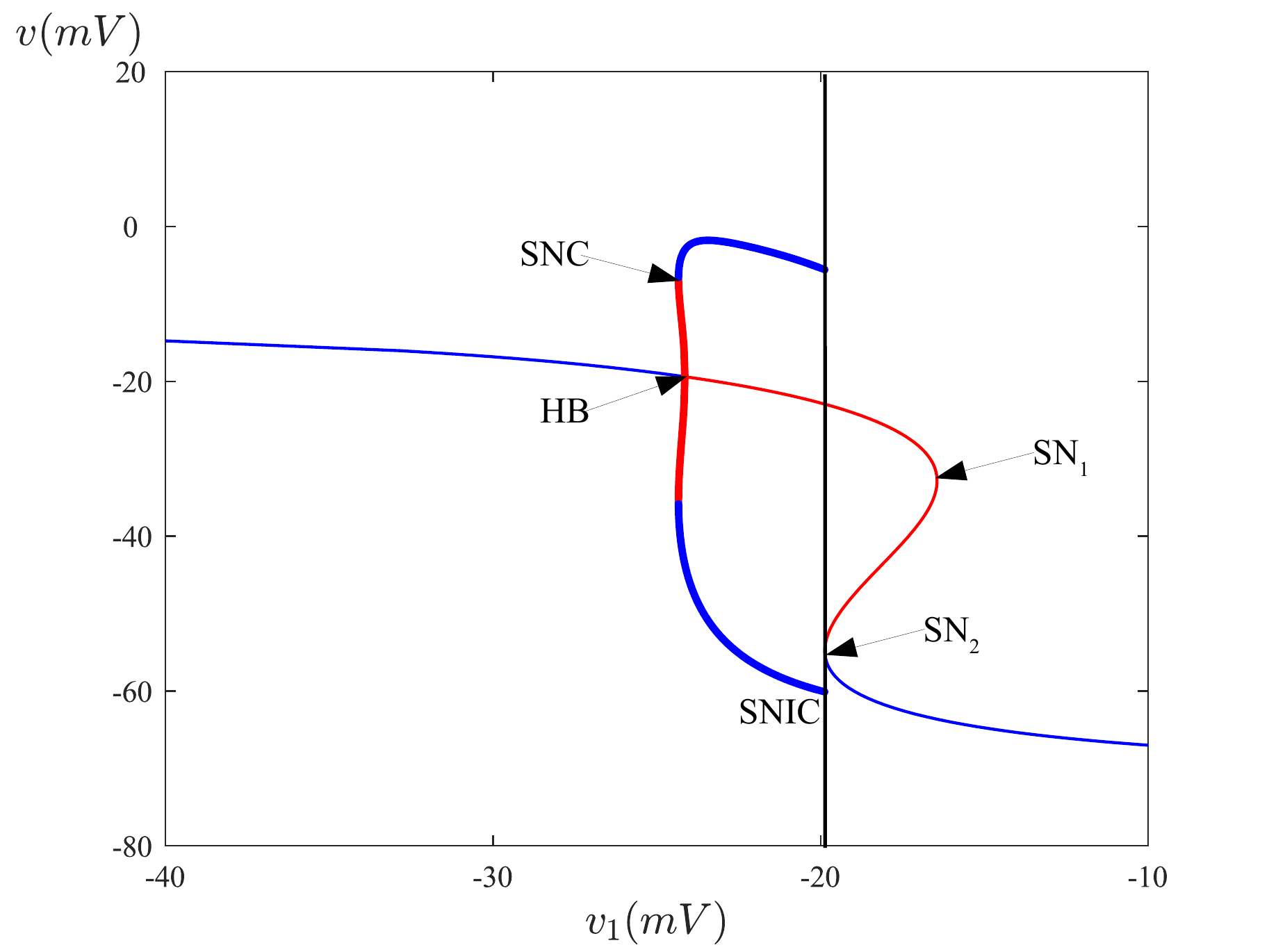}
  \label{fig:2.4b}
\end{subfigure}\\%
\begin{subfigure}{.5\textwidth}
  \centering
  \caption{}
  \includegraphics[width=\textwidth]{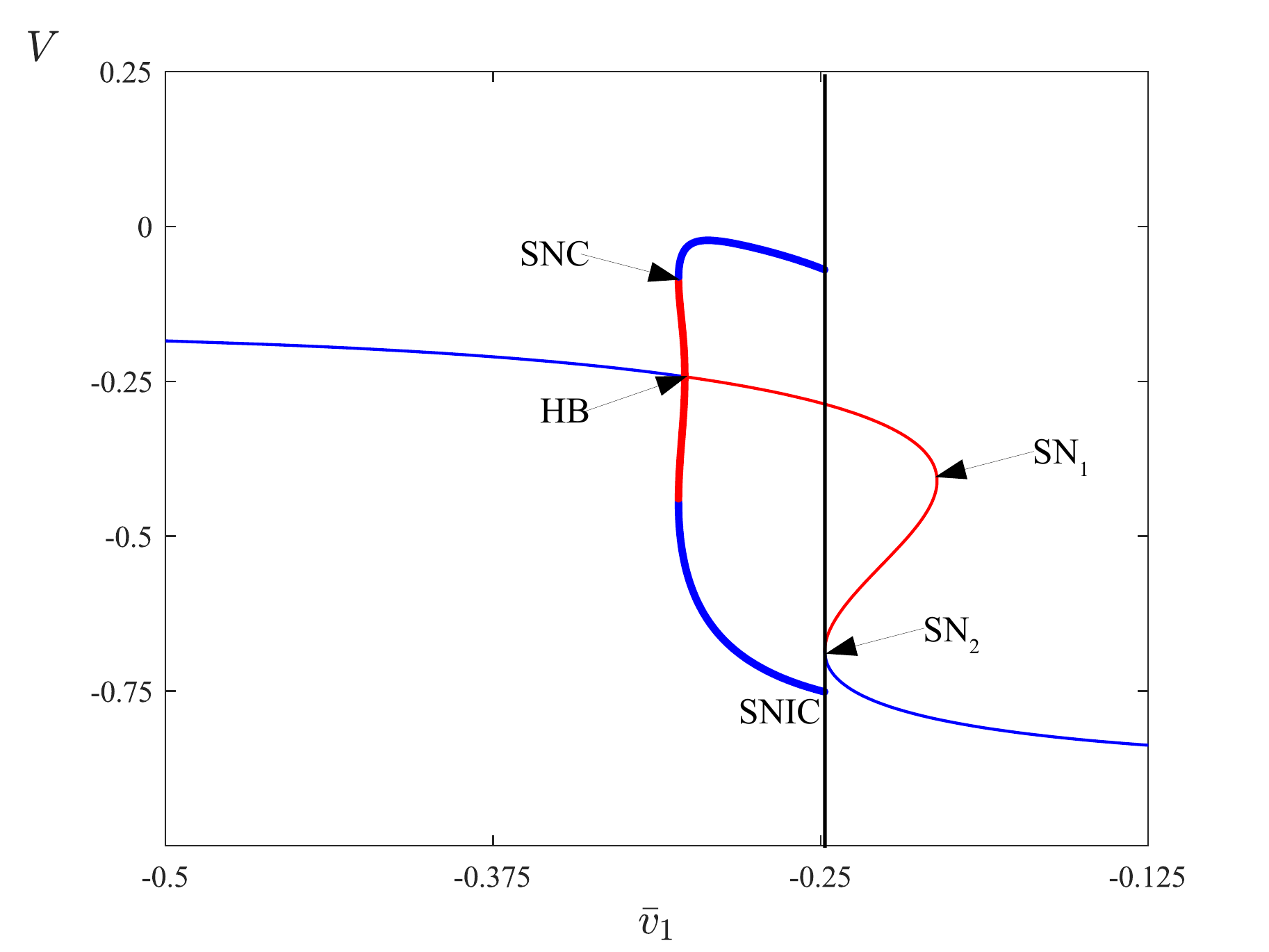}
  \label{fig:2.4c}
\end{subfigure}%
\begin{subfigure}{.5\textwidth}
      \centering
      \caption{}
      \includegraphics[width=\textwidth]{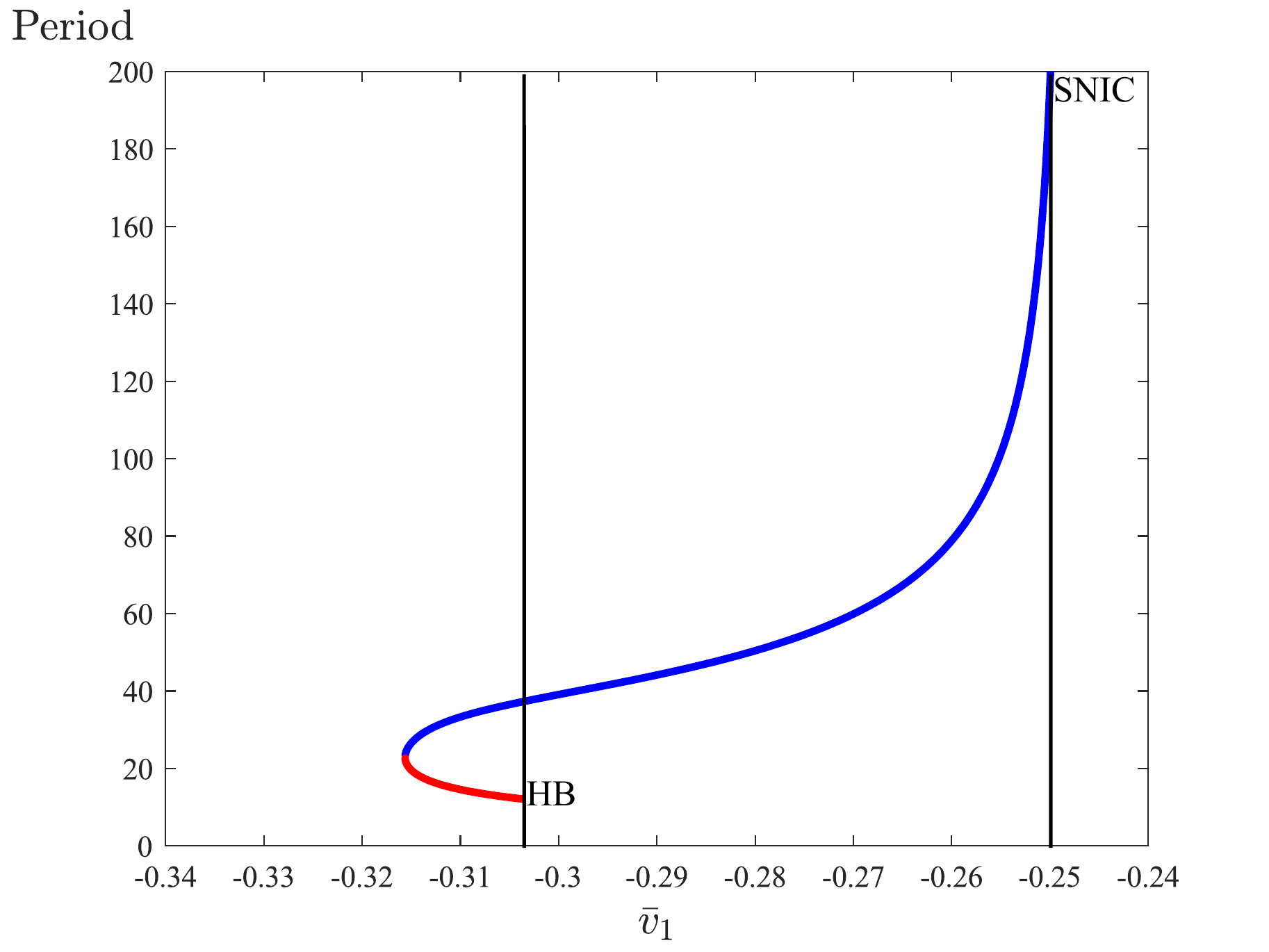}
      \label{fig:2.4d}
\end{subfigure}
\caption{Bifurcation diagrams of (a) the
full model \eqref{eq:V}--\eqref{eq:Ca} with $v_1$ as the bifurcation parameter, (b) the reduced model \eqref{eq:reduced_2}--\eqref{eq:reduced_3}
with $v_1$ as the bifurcation parameter and (c) the nondimensionalised model \eqref{eq:dimless1}--\eqref{eq:dimless2} with $\bar{v}_1$ as the bifurcation parameter.
The remaining parameter values are given in Tables \ref{Tab:1} and \ref{Tab:2}. Panel (d) shows the period of the oscillations for the nondimensionalised model. The blue and red curves represent stable  and unstable periodic orbits in Fig.~\ref{fig:2.4c}.
Thin [thick] curves correspond to equilibria [periodic orbits].
Blue [red] curves correspond to stable [unstable] solutions.
HB: Hopf bifurcation;
SN: saddle-node bifurcation (of an equilibrium); SNC: saddle-node bifurcation of a periodic orbit;
SNIC: saddle-node on an invariant circle bifurcation.}
\label{fig:4}
\end{figure}

Next, we verify the excitability property of the model by varying the voltage associated with the fraction of open $\text{K}^{+}$ channels as a bifurcation parameter. Since  $v_1$ is dependent on transmural pressure \citep{Gonzalez-Fernandez1994}, it is considered to be the main bifurcation parameter in the full model. For the reduced model this parameter is $\bar{v}_1$. We choose a range of values of $v_{1}$ and $\bar{v}_1$ for which the systems either converge to a steady state (absence of vasomotion) or oscillate (presence of vasomotion). We use values of $v_1$ between $-40$mV and $-10$mV, which corresponds to values of $\bar{v}_1$ between $-0.5$ and $-0.125$. Figs.~\ref{fig:2.4a}--\ref{fig:2.4c} show the bifurcation diagrams of the full, reduced and nondimensionalised models. A detailed discussion of the bifurcation diagrams, particularly for the nondimensionalised model is given in Sect.~\ref{sec:bifurcation}. 

 \section{Bifurcation analysis of Type I and Type II excitability}\label{sec:bifurcation}

Here we investigate the dynamics of the nondimensionalised model \eqref{eq:dimless1}--\eqref{eq:dimless2} via a bifurcation analysis. In Sect.~\ref{sec:codimension-1} the influence of different model parameters on model behaviour is considered. Then in Sect.~\ref{sec:codimesion-2} we relate transitions between Type I and Type II excitability to codimension-two bifurcations.

\subsection{Changes to the dynamics as one parameter is varied}
\label{sec:codimension-1}
As shown in Fig.~\ref{fig:2.3c}
the nondimensionalised model exhibits stable oscillations for the parameter values of Table \ref{Tab:2}.
Here we study how the dynamics changes as the parameters
$\bar{v}_1$, $\bar{v}_3$, and $\bar{v}_L$ are varied from their values in Table \ref{Tab:2}. First we consider $\bar{v}_1$.
A bifurcation diagram is shown in Fig.~\ref{fig:2.4c}.
We observe the system has a unique equilibrium except between two saddle-node bifurcations,
$\text{SN}_1$ and $\text{SN}_2$.
To the right of $\text{SN}_2$ the lower equilibrium branch is the only stable solution of the system.
The saddle-node bifurcation $\text{SN}_2$ is in fact a SNIC bifurcation
(saddle-node on an invariant circle)
as here there exists an orbit homoclinic to the equilibrium
\citep{KuznetsovY.A.1995ElementsTheory}
.
To the left of $\text{SN}_2$ this orbit persists as a stable periodic orbit. Thus here \eqref{eq:dimless1}--\eqref{eq:dimless2} model SMC activity with Type I excitability
\citep{Hodgkin1952, Ermentout1996TypeSychrony, Izhikevich2007DynamicalBursting}.
As we pass through the SNIC bifurcation by decreasing the value of $\bar{v}_1$
the excitable state changes to periodic oscillations.
As shown in Fig.~\ref{fig:2.4d}
the period of the oscillations decreases from infinity
as a consequence of the homoclinic connection.\\

Upon further decrease in the value of $\bar{v}_1$
the stable periodic orbit loses stability in a saddle-node bifurcation (SNC).
The resulting branch of unstable periodic orbits
terminates in a subcritical Hopf bifurcation (HB).
Between these bifurcations the system is bistable because
the upper equilibrium branch is stable to the left of the Hopf bifurcation.\\

Next we vary the value of the parameter $\bar{v}_3$.  This is because it is of biological interest to understand the influence of transmural pressure. In the full model \eqref{eq:V}--\eqref{eq:Ca} transmural pressure is associated with the parameter $v_6$, so in the nondimensionalised model it is associated with $\bar{v}_3$ through $v_3^* = v_6 + \frac{v_5}{2}$. Hence we can examine the influence of transmural pressure by using $\bar{v}_3$ as a bifurcation parameter.\\

\begin{figure}
\centering
\begin{subfigure}{.5\textwidth}
  \centering
  \caption{}
  \includegraphics[width = \textwidth]{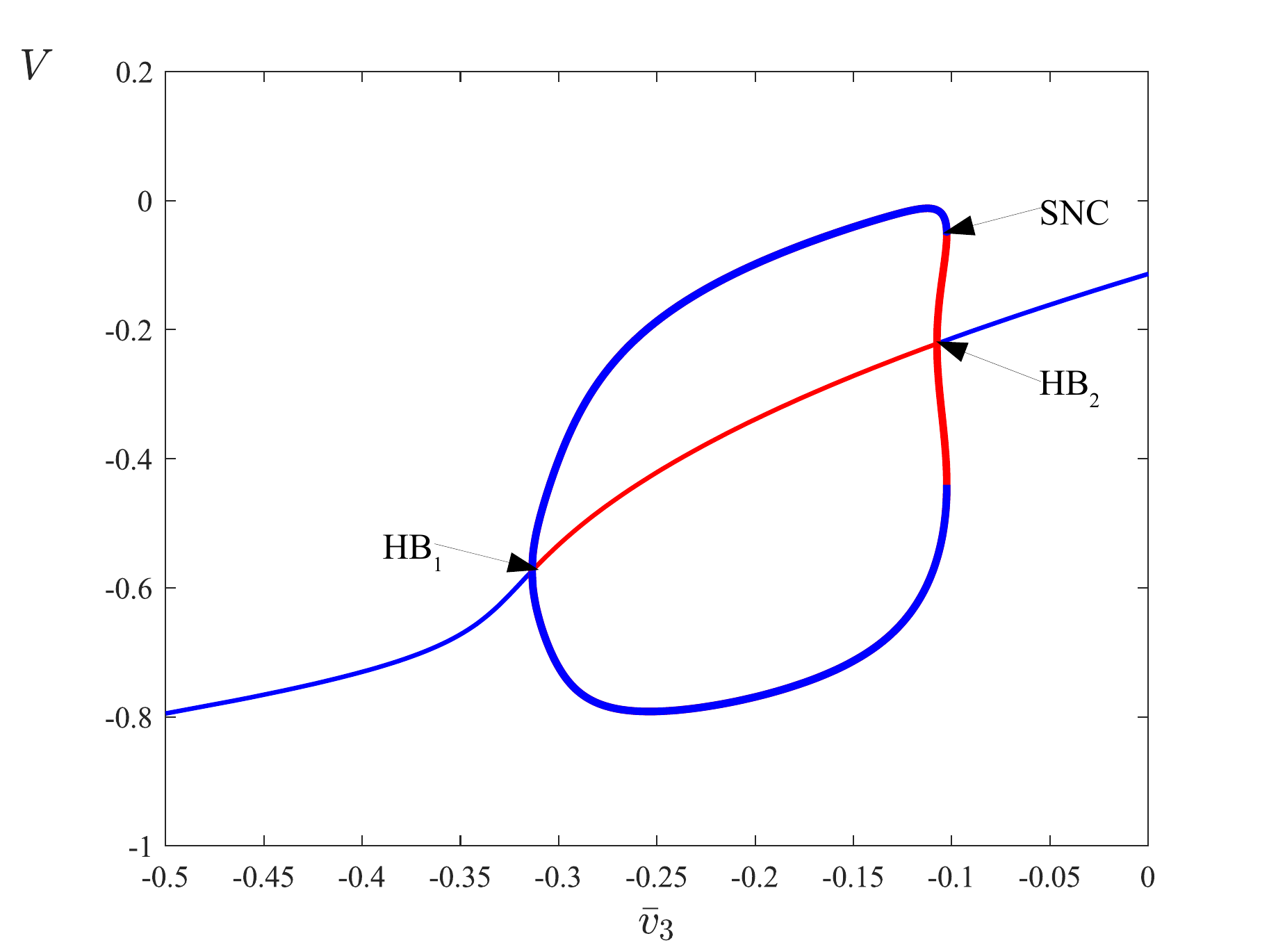}
  \label{fig:5a}
\end{subfigure}%
\begin{subfigure}{.5\textwidth}
  \centering
  \caption{}
  \includegraphics[width = \textwidth]{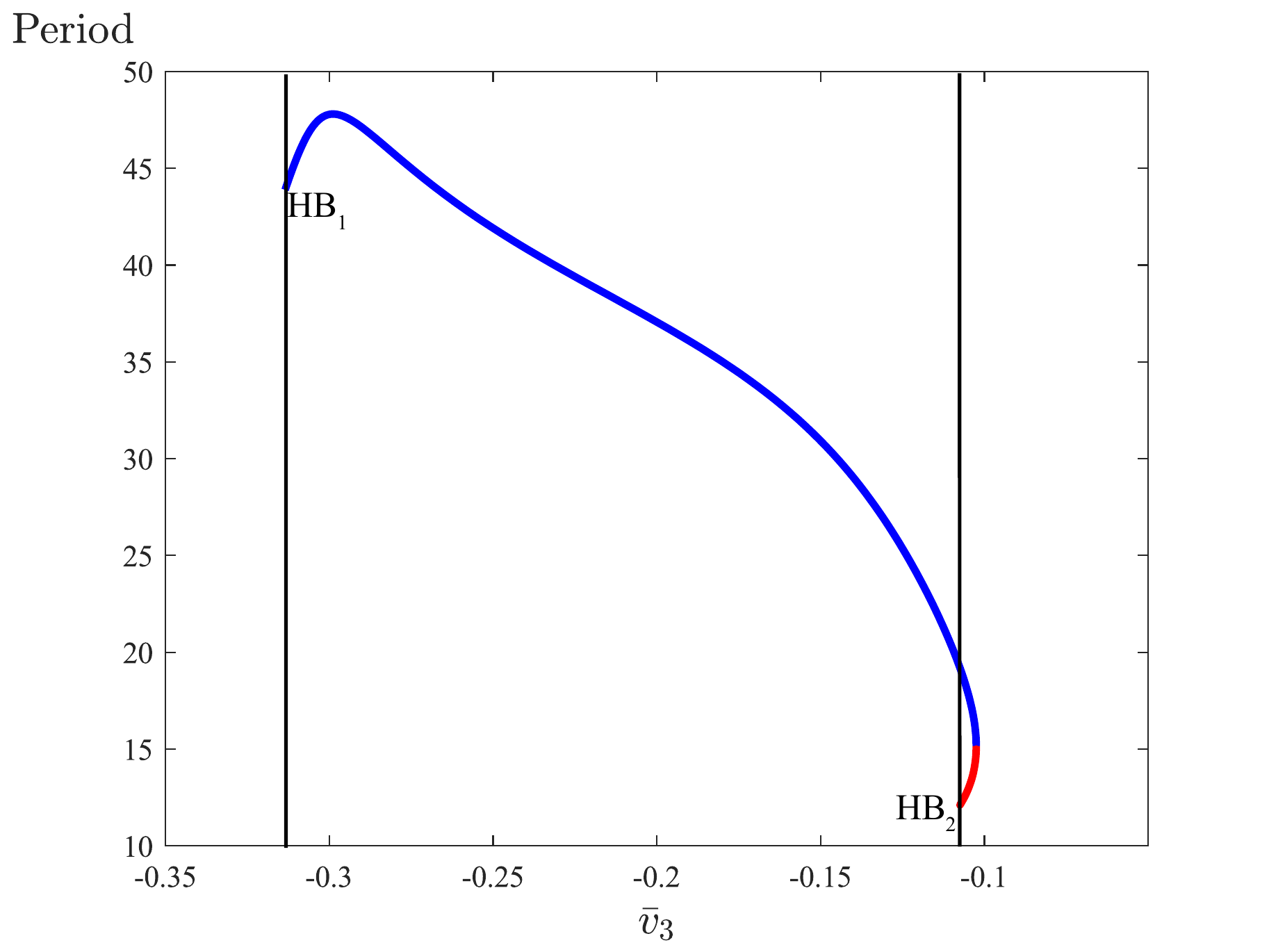}
  \label{fig:5b}
\end{subfigure}

\caption{(a) A bifurcation diagram of the nondimensionlised model  \eqref{eq:dimless1}--\eqref{eq:dimless2} with $\bar{v}_{3}$ as the bifurcation parameter
and other parameter values as given in Table \ref{Tab:2}.
(b) A plot of the periodic oscillations as a function of parameter $\bar{v}_3$.
The labels and other conventions are as in Fig.~\ref{fig:4}.
}
\label{fig:5}
\end{figure}

As shown in Fig.~\ref{fig:5a}, as we increase the value of $\bar{v}_3$ a unique equilibrium loses stability in a supercritical Hopf bifurcation $\text{HB}_1$ then regains stability in a subcritical Hopf bifurcation $\text{HB}_2$.  Therefore in this case the system exhibits Type II excitability. The stable oscillations are created in $\text{HB}_1$ with finite period (see Fig.~\ref{fig:5b}).
They subsequently lose stability at the saddle-node bifurcation SNC
and terminate at $\text{HB}_2$.

\begin{figure}
\centering
  \includegraphics[width=0.6\textwidth]{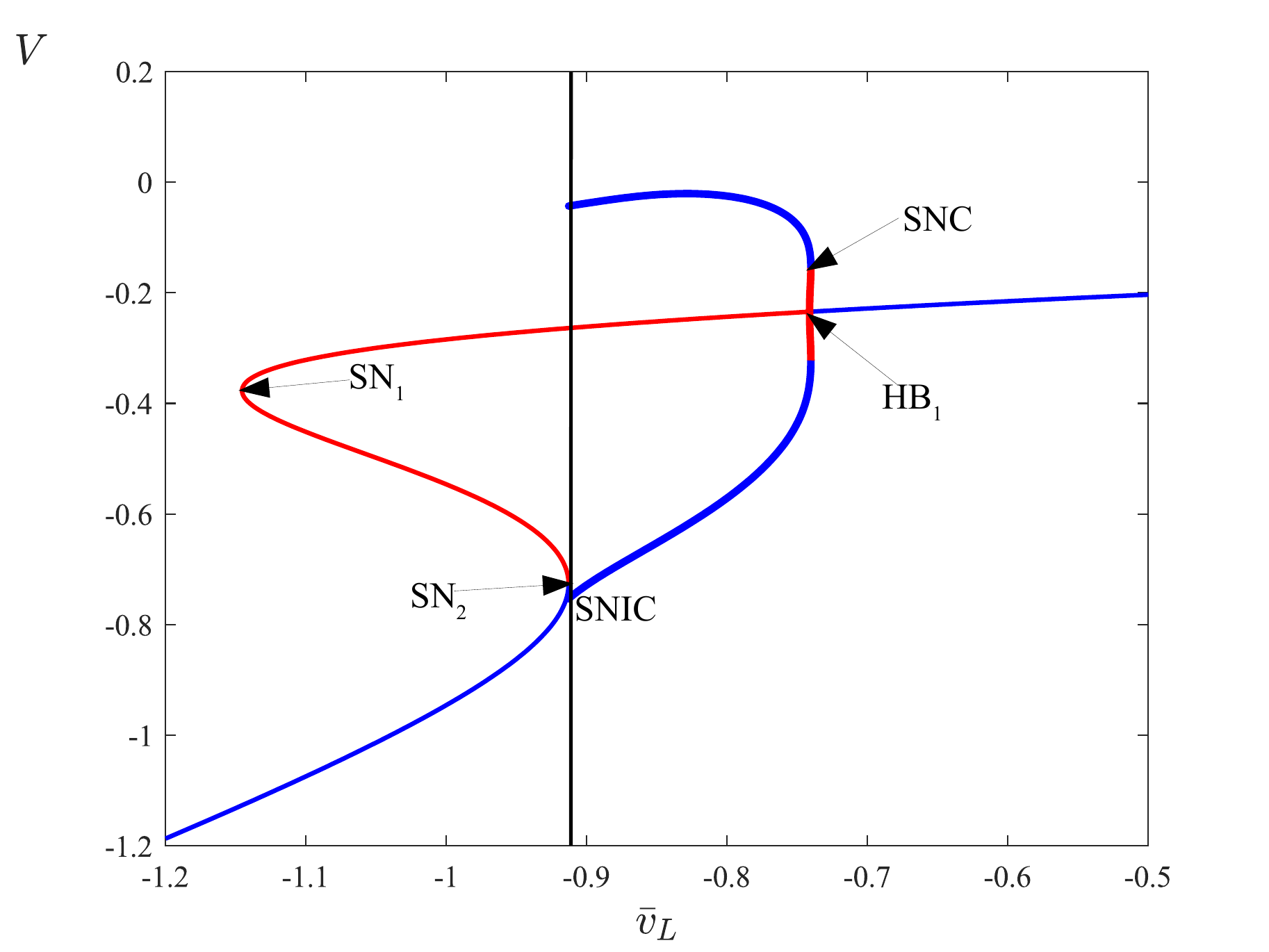}
\caption{A bifurcation diagram of the nondimensionalised model \eqref{eq:dimless1}--\eqref{eq:dimless2} with $\bar{v}_L$ as the bifurcation parameter 
and other parameter values as given in Table \ref{Tab:2}.
The labels and other conventions are as in Fig.~\ref{fig:4}}
  \label{fig:6}
\end{figure}

Lastly, variation of $\bar{v}_L$ produces the bifurcation diagram Fig.~\ref{fig:6}.
This has the same type of bifurcation structure as Fig.~\ref{fig:2.4b} (except in reverse).
Thus increasing the value of $\bar{v}_L$ results in the same qualitative changes to the dynamics
as decreasing the value of $\bar{v}_1$. In particular the excitability is Type I.

\subsection{Transitions between types of excitability}\label{sec:codimesion-2}
\begin{figure}
 \centering
  \includegraphics[width =0.9 \textwidth]{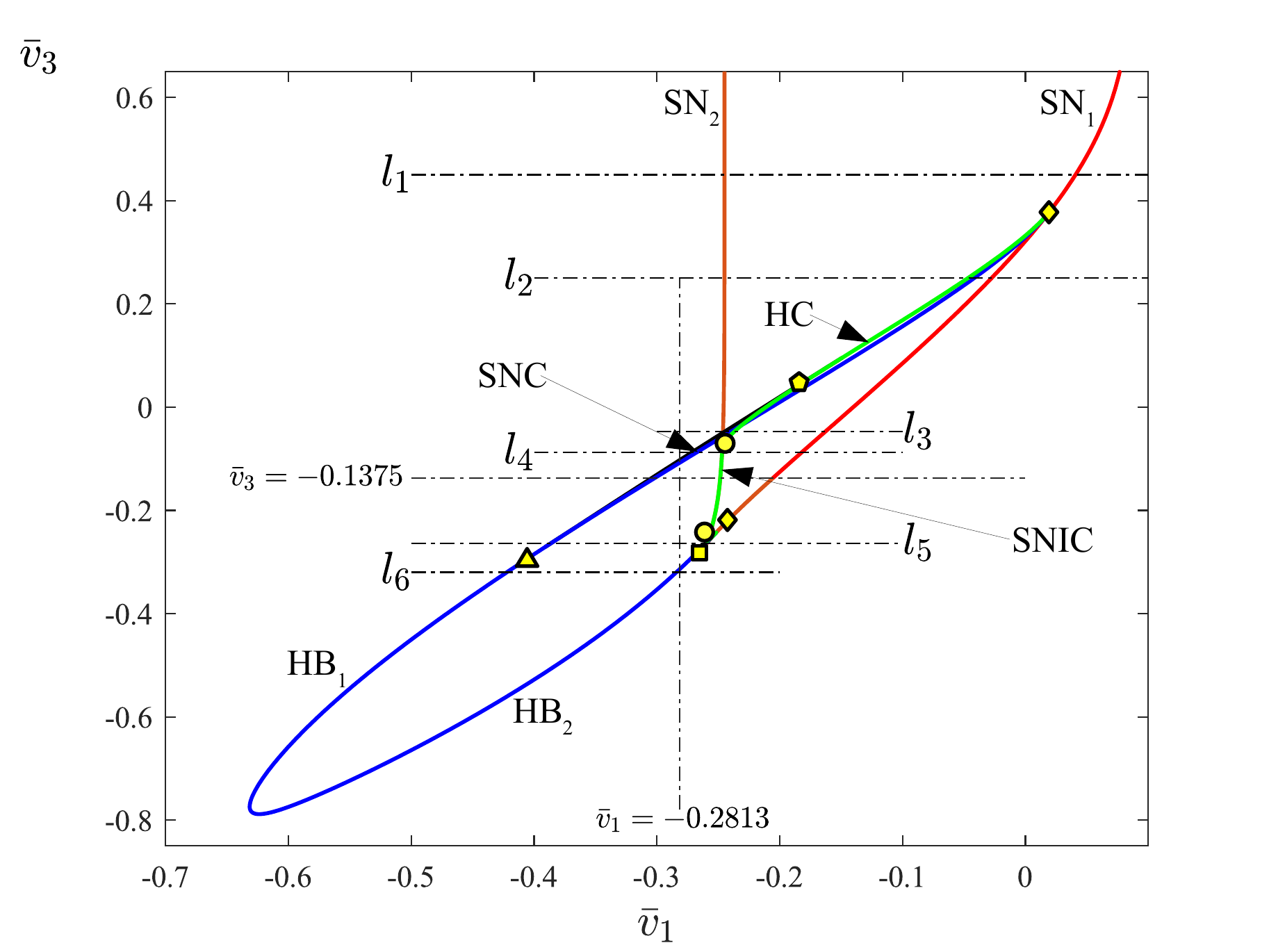}
   \caption{ A two-parameter bifurcation diagram of the nondimensionalised model \eqref{eq:dimless1}--\eqref{eq:dimless2} in the $(\bar{v}_1,\bar{v}_3)$-plane for the parameter values of Table \ref{Tab:2}. The values of $\bar{v}_3$ in $l_1$, $l_2$, $l_3$, $l_4$, $l_5$ and $l_6$ are $0.45$, $0.25$, $-0.047$, $-0.088$, $-0.26$ and $-0.32$, respectively. The loci of codimension-one bifurcations are coloured as follows: blue: Hopf bifurcation, red: saddle-node bifurcation (or SNIC), green: homoclinic bifurcation, and black: saddle-node bifucation of periodic orbit. The labels for the codimension-two bifurcations are explained in Table \ref{Tab:3}}
  \label{fig:7}
\end{figure}

In this section we perform a two-parameter bifurcation analysis of the nondimensionalised model
\eqref{eq:dimless1}--\eqref{eq:dimless2} by varying the parameters $\bar{v}_1$ and $\bar{v}_3$.
This is summarised by the two-parameter bifurcation diagram, Fig.~\ref{fig:7},
which was produced via the numerical continuation software AUTO-07p \citep{Doedel2012}.
Two of the one-parameter bifurcation diagrams described above,
are slices of Fig.~\ref{fig:7}.
Specifically Fig.~\ref{fig:2.4c} has the value of $\bar{v}_3$ fixed at $-0.1375$
and Fig.~\ref{fig:5a} has the value of $\bar{v}_1$ fixed at $-0.2813$.

In the remainder of this section we describe Fig.~\ref{fig:7}
and consequences to transitions between Type I and II excitability
by studying slices at six different values of $\bar{v}_3$.
Fig.~\ref{fig:7} includes five different codimension-two bifurcations
summarised by Table.~\ref{Tab:3} and discussed below.

\def\brcurs{{\mbox{$\includegraphics[scale=0.3]{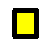} $}}}
\def\prcurs{{\mbox{$\includegraphics[scale=0.3]{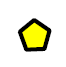} $}}}
\def\hrcurs{{\mbox{$\includegraphics[scale=0.3]{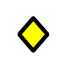} $}}}
\def\trcurs{{\mbox{$\includegraphics[scale=0.3]{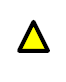} $}}}
\def\crcurs{{\mbox{$\includegraphics[scale=0.3]{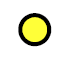} $}}}

\begin{table}
\centering
\caption{Abbreviations and notations  of codimension-two bifurcations}
\label{Tab:3}       
\begin{tabular}{lll}
\hline\noalign{\smallskip}
Bifurcation & Abbreviation & Label \\
\noalign{\smallskip}\hline\noalign{\smallskip}
Cusp bifurcation & CP & \hspace{3.5mm} \raisebox{-1mm}[0mm][0mm]{\brcurs} \\[1mm]
Bogdanov-Takens bifurcation & $\text{BT}_i$ \hspace{0.5mm}$i=1,2$ & \hspace{2.85mm} \raisebox{-1.1mm}[0mm][0mm]{\hrcurs} \\[1mm]
Generalized Hopf bifurcation & GH & \hspace{2.8mm} \raisebox{-2.3mm}[0mm][0mm]\trcurs\\[1mm]
Resonant homoclinic bifurcation & RHom & \hspace{3.5mm}\raisebox{-1.2mm}[0mm][0mm]\prcurs  \\[1mm]
Non-central saddle-node homoclinic bifurcation & $\text{NSH}_i$ \hspace{0.5mm}$i=1,2$ &\hspace{3.5mm}\raisebox{-1.3mm}[0mm][0mm] \crcurs\\[1mm]
\noalign{\smallskip}\hline
\end{tabular}
\end{table}
 
\begin{figure}
\centering
\begin{subfigure}{.6\textwidth}
 \centering
  \caption{}
  \includegraphics[width = \textwidth]{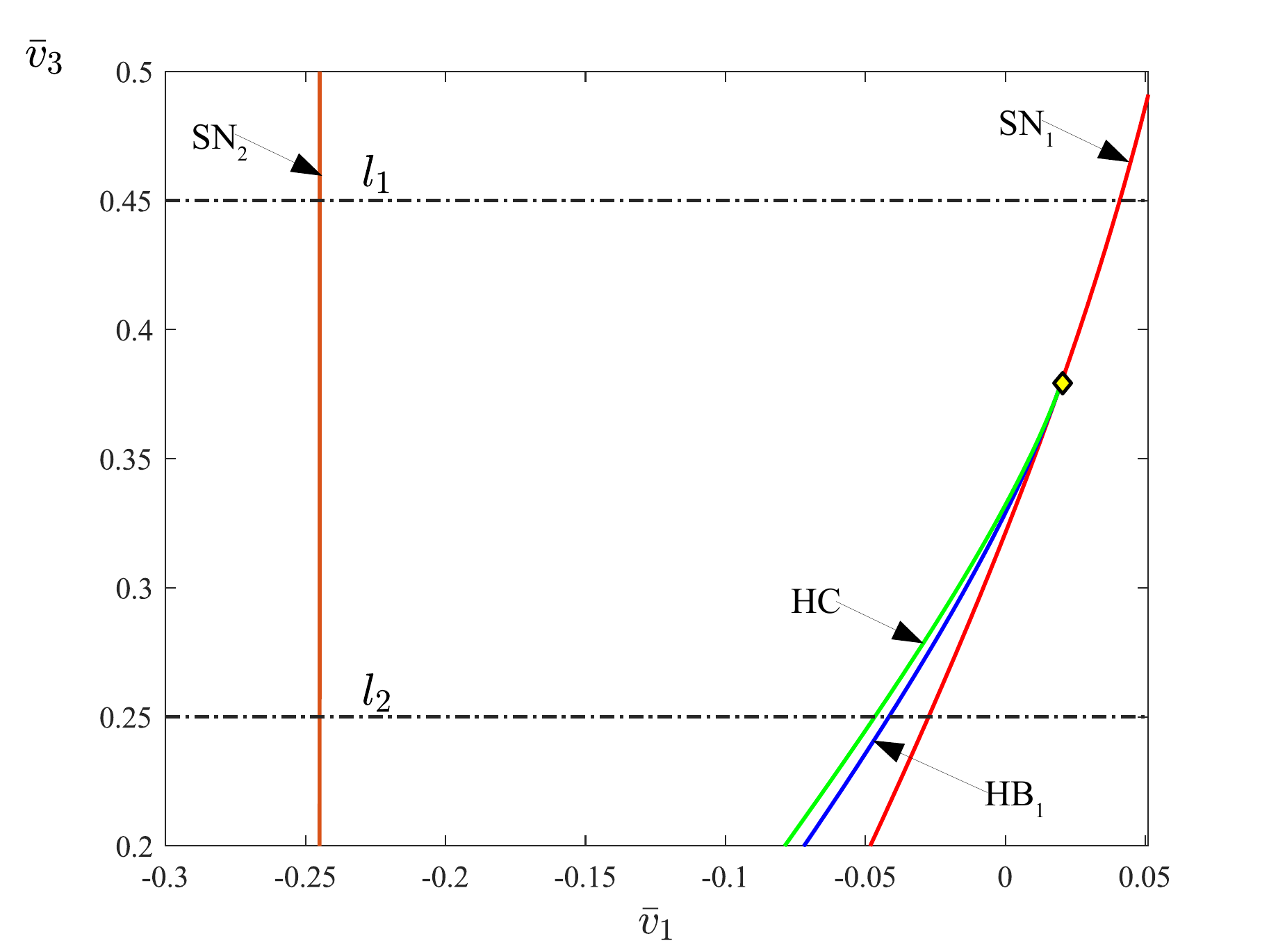}
  \label{fig:8a}
\end{subfigure}
\begin{subfigure}{.5\textwidth}
 \centering
  \caption{}
  \includegraphics[width = \textwidth]{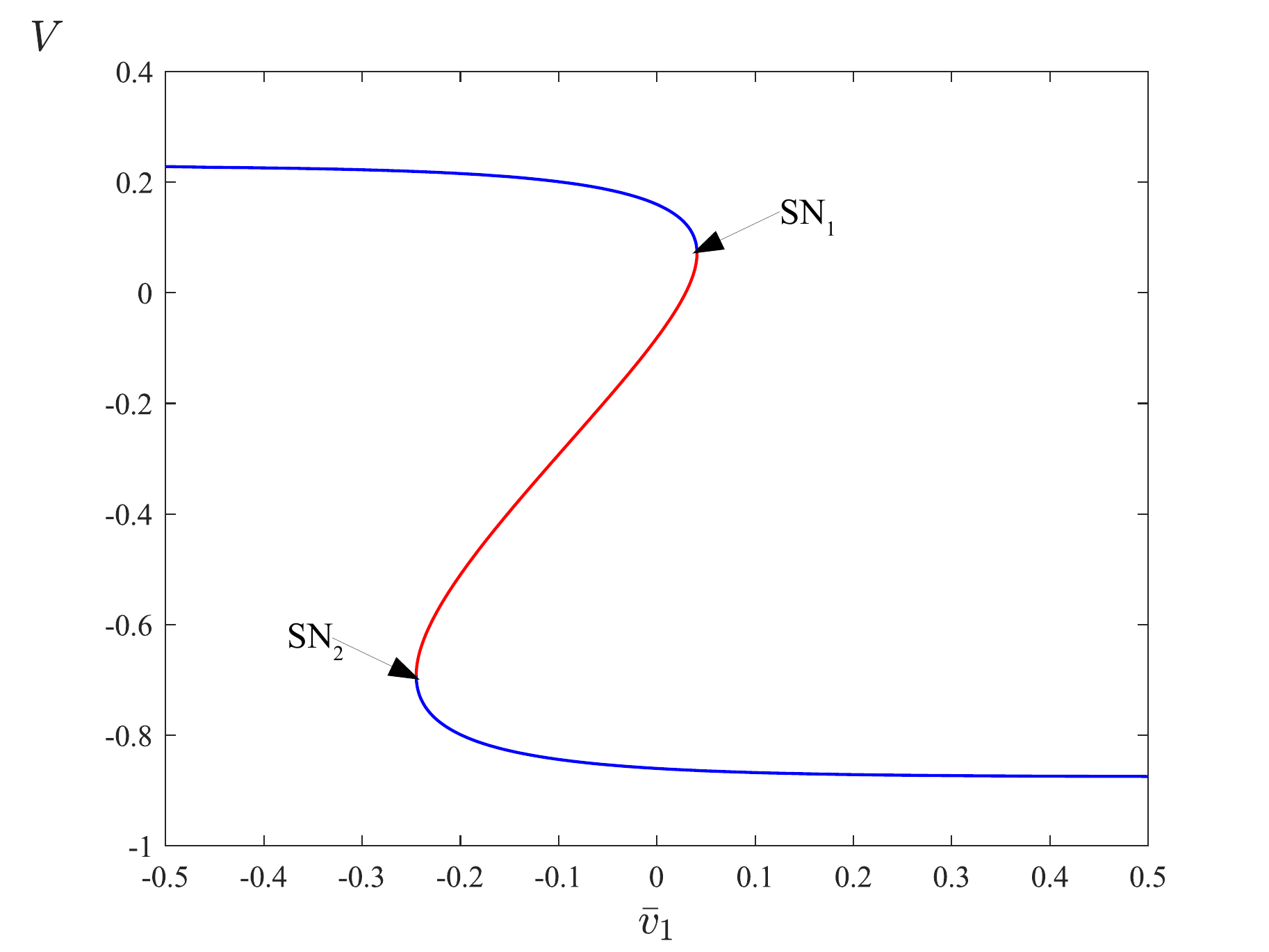}
  \label{fig:8b}
\end{subfigure}%
\begin{subfigure}{.5\textwidth}
  \centering
  \caption{}
  \includegraphics[width = \textwidth]{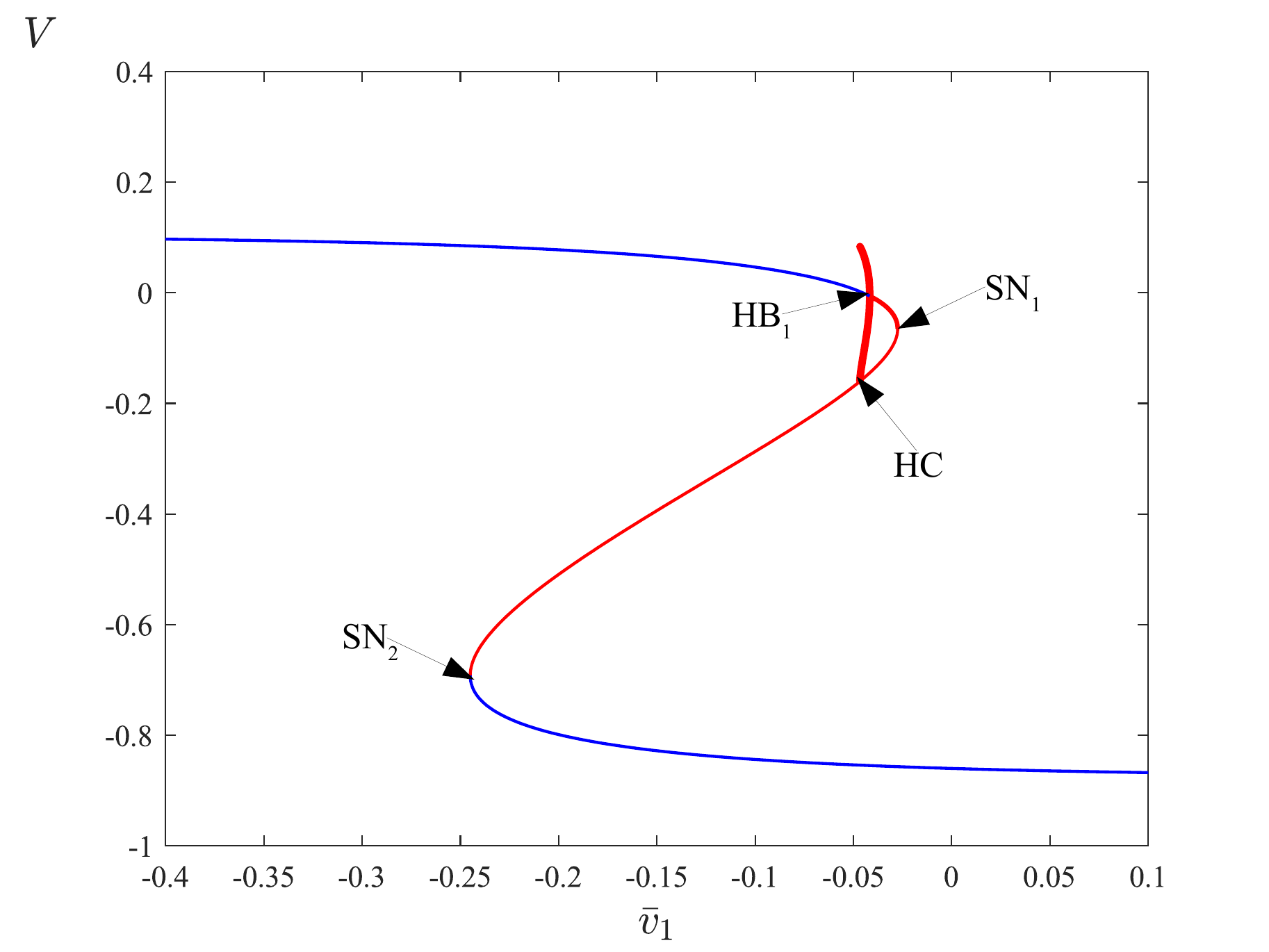}
  \label{fig:8c}
\end{subfigure}
\caption{(a) An enlargement of Fig. \ref{fig:7} showing lines $l_1$ and $l_2$. The filled diamond is a Bogdanov-Takens bifurcation.. (b) A one-parameter bifurcation diagram  along $l_1$ with $\bar{v}_3=0.45$. (c) A one-parameter bifurcation diagram along $l_2$ with $\bar{v}_3=0.25$. HB: Hopf bifurcation, $\text{SN}$: saddle-node bifurcation, SNC: saddle-node bifurcation of a periodic orbit, HC: homoclinic bifurcation}
\label{fig:8}
\end{figure}

For sufficiently large values of $\bar{v}_3$ the only bifurcations are the two saddle-node
bifurcations $\text{SN}_1$ and $\text{SN}_2$, 
see Fig.~\ref{fig:8a} which shows a magnification of Fig.~\ref{fig:7}.
Thus for the slice $l_1$ there are no periodic solutions, Fig.~\ref{fig:8b}

As we decrease the value of $\bar{v}_3$ a Bogdanov-Takens bifurcation
\citep{Takens1974SINGULARITIESFIELDS, Bogdanov1975VersalEigenvalues}, denoted $\text{BT}_1$, occurs
on the saddle-node locus $\text{SN}_1$ at $\bar{v}_3 \approx 0.3792$.
This is a codimension-two point from which loci
of homoclinic and subcritical Hopf bifurcations emanate, denoted $\text{HC}$ and $\text{HB}_1$.
As known from the theory of Bogdanov-Takens bifurcations
\citep{KuznetsovY.A.1995ElementsTheory}
and as seen in Fig.~\ref{fig:8a}
these loci are tangent to $\text{SN}_1$ at the codimension-two point.
Thus for a slice below $\text{BT}_1$, such as $l_2$ for which $\bar{v}_3 = 0.25$,
apart from the saddle-node bifurcations already observed
there are now also homoclinic and Hopf bifurcations
between which there exists an unstable periodic orbit, Fig.~\ref{fig:8c}.
Observe also that upon crossing $\text{BT}_1$ the interval of values of $\bar{v}_1$
in which the system is bistable changes from endpoints at $\text{SN}_2$ and $\text{SN}_1$ (for $l_1$)
to endpoints at $\text{SN}_2$ and $\text{HB}_1$ (for $l_2$).

\begin{figure}
\centering
  \begin{subfigure}[b]{.6\textwidth}
    \centering
    \caption{}
    \includegraphics[width = \textwidth]{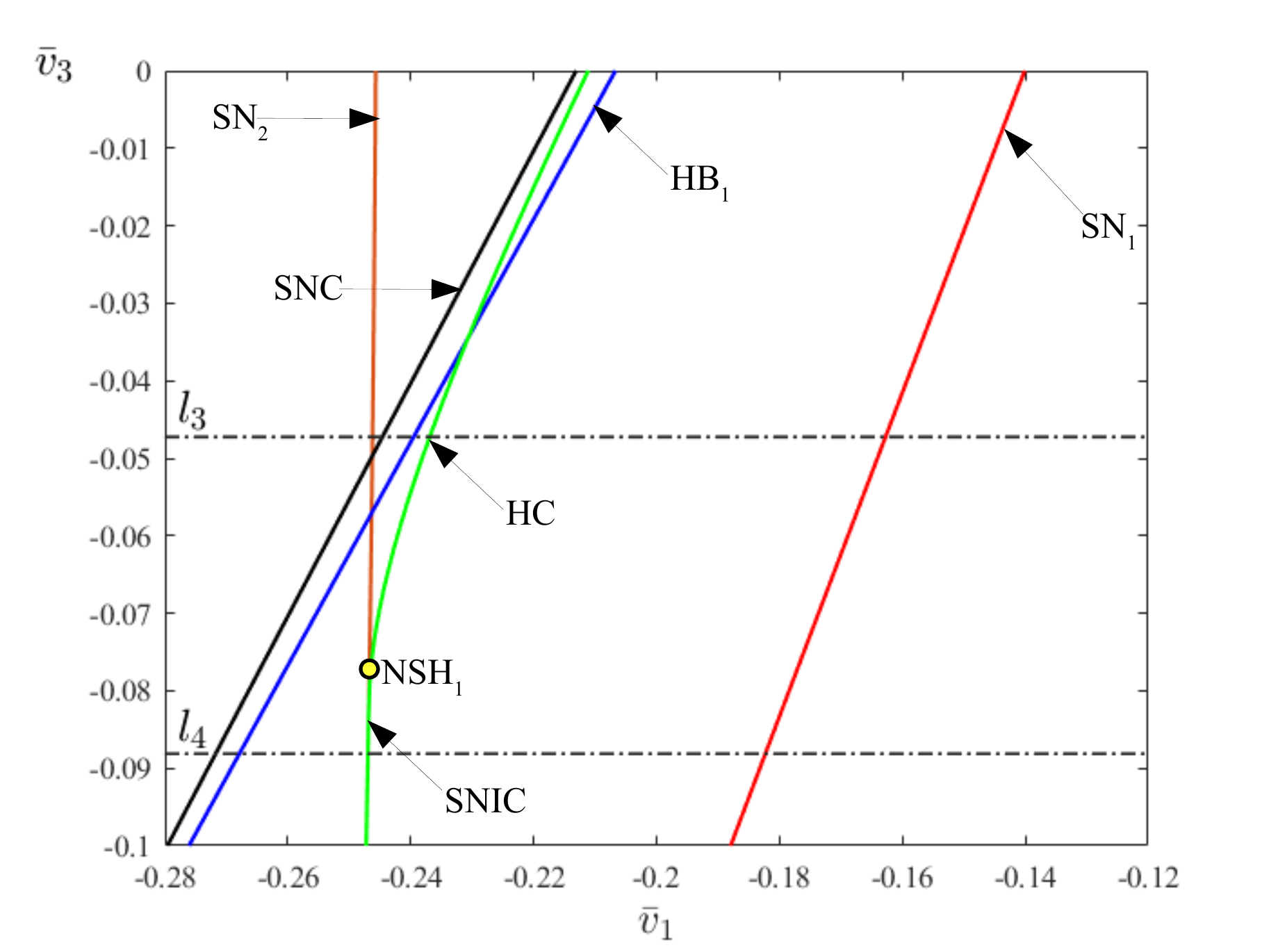}
 \label{fig:9a}
  \end{subfigure}\\%
  \begin{subfigure}[b]{.5\textwidth}
    \centering
    \caption{}
    \includegraphics[width = \textwidth]{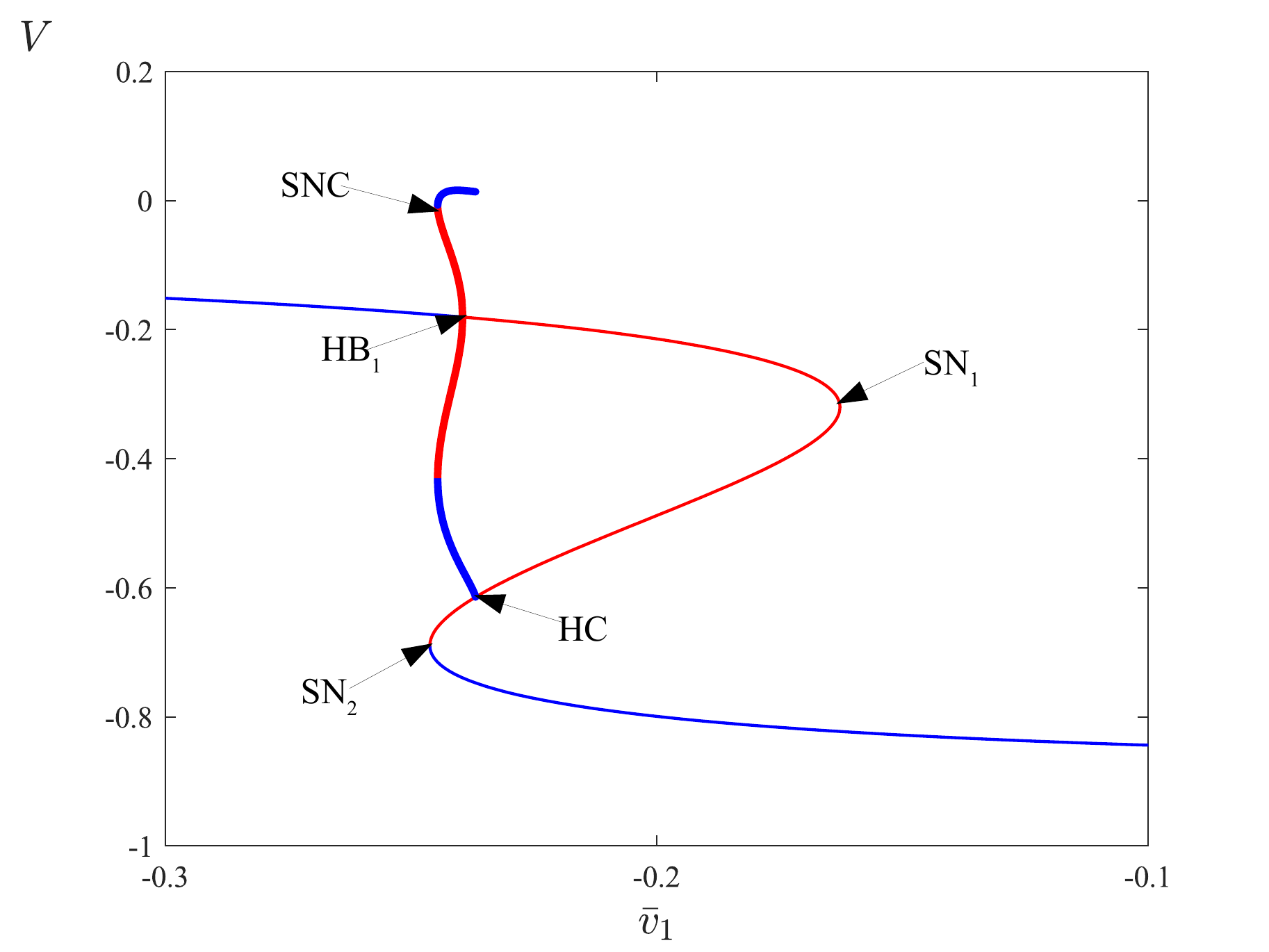}
    \label{fig:9b}
  \end{subfigure}%
  \begin{subfigure}[b]{.5\textwidth}
    \centering
    \caption{}
    \includegraphics[width = \textwidth]{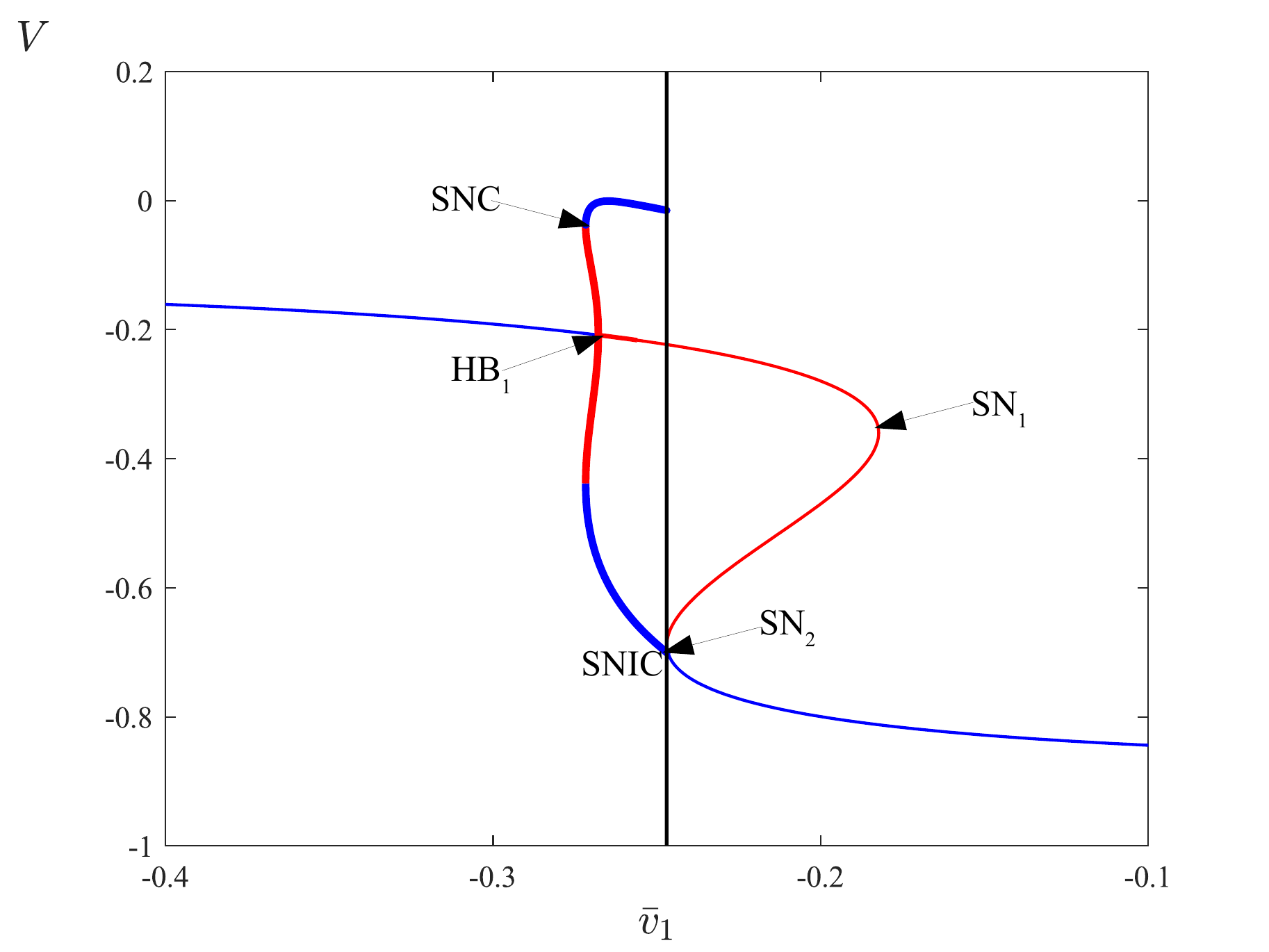}
    \label{fig:9c}
  \end{subfigure}%
 \caption{ (a) An enlargement of Fig. \ref{fig:7}  showing lines $l_3$ and $l_4$. The filled circle is a non-central saddle-node homoclinic bifurcation. (b) A one-parameter bifurcation diagram along $l_3$ with $\bar{v}_3=-0.047$. (c) A one-parameter bifurcation diagram along  $l_4$ with $\bar{v}_3=-0.088$. HB: Hopf bifurcation, $\text{SN}$: saddle-node bifurcation, SNC: saddle-node bifurcation of a periodic orbit, SNIC: saddle-node on an invariant circle bifurcation, HC: homoclinic bifurcation}
\label{fig:9}
\end{figure}

\begin{figure}
 \centering
  \includegraphics[width = 0.7\textwidth]{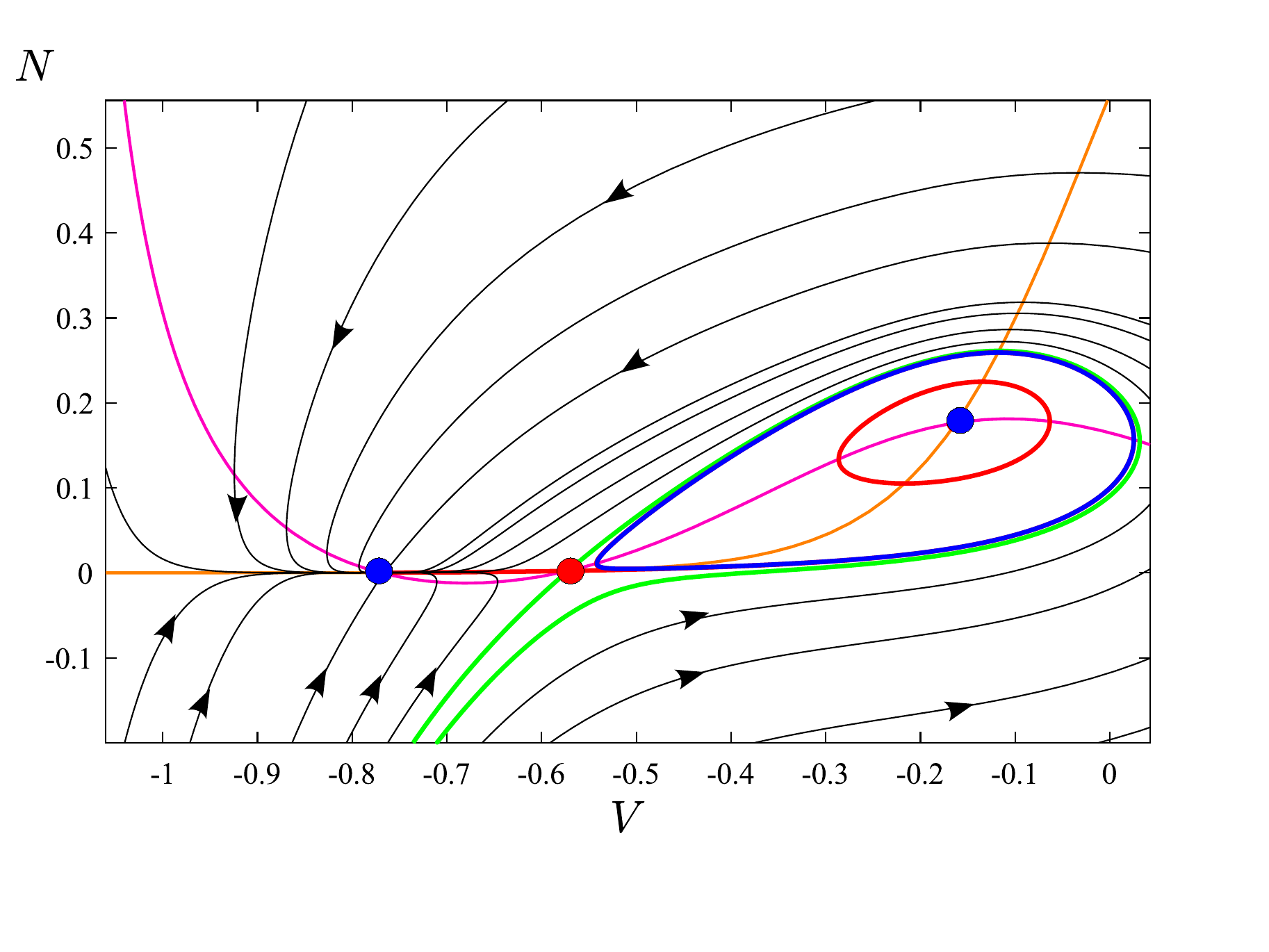}
   \caption{ A phase portrait of the nondimensionalised model \eqref{eq:dimless1}--\eqref{eq:dimless2} on line $l_3$ at $\bar{v}_3 =-0.047$ showing tristability. The blue and red curves are stable and unstable periodic orbits. The magenta and orange curves are the nullclines for $N$ and $V$. The black curves are the solution tractories. The blue and red circles are stable and unstable equilibria}
  \label{fig:10}
\end{figure}

As the value of $\bar{v}_3$ is decreased further, $\text{HB}_1$ shifts to the left
and a locus of saddle-node bifurcations of the periodic orbit, $\text{SNC}$,
emanates from the codimension-two point $\text{RHom}$
on $\text{HC}$ at $\bar{v}_3 \approx 0.0095$, see Fig.~\ref{fig:9a}.
Thus below this point there exists a stable periodic orbit
between $\text{SNC}$ and $\text{HC}$, such as for the slice $l_3$, Fig.~\ref{fig:9b}.
For this slice, as the value of $\bar{v}_1$ is decreased
stable oscillations are created at $\text{HC}$.
Here there is a small region of tristability:
stable oscillations coexist with two stable equilibria, see Fig.~\ref{fig:10}.

Upon further decrease of $\bar{v}_3$ the locus $\text{HC}$ collides tangentially with $\text{SN}_2$
at the codimension-two point $\text{NSH}_1$.
This is known as a non-central saddle-node homoclinic bifurcation,
see for instance \citep{Govaerts2005TheApproach}.
The collision produces the locus $\text{SNIC}$ (saddle-node of an invariant circle).
Thus immediately below $\text{NSH}_1$ the system exhibits Type I excitability.
The system transitions from a stable equilibrium to a stable periodic orbit
at the SNIC bifurcation, such as for the slice $l_4$, Fig.~\ref{fig:9c}
(and as described earlier, Fig.~\ref{fig:2.4b}).
Thus the point $\text{NSH}_1$ marks the onset of Type I excitability.
This has been observed previously for the reduced Morris-Lecar model with external current
\citep{Tsumoto2006BifurcationsModel}.

\begin{figure}
\centering
  \begin{subfigure}[b]{.56\textwidth}
    \centering
    \caption{}
    \includegraphics[width = \textwidth]{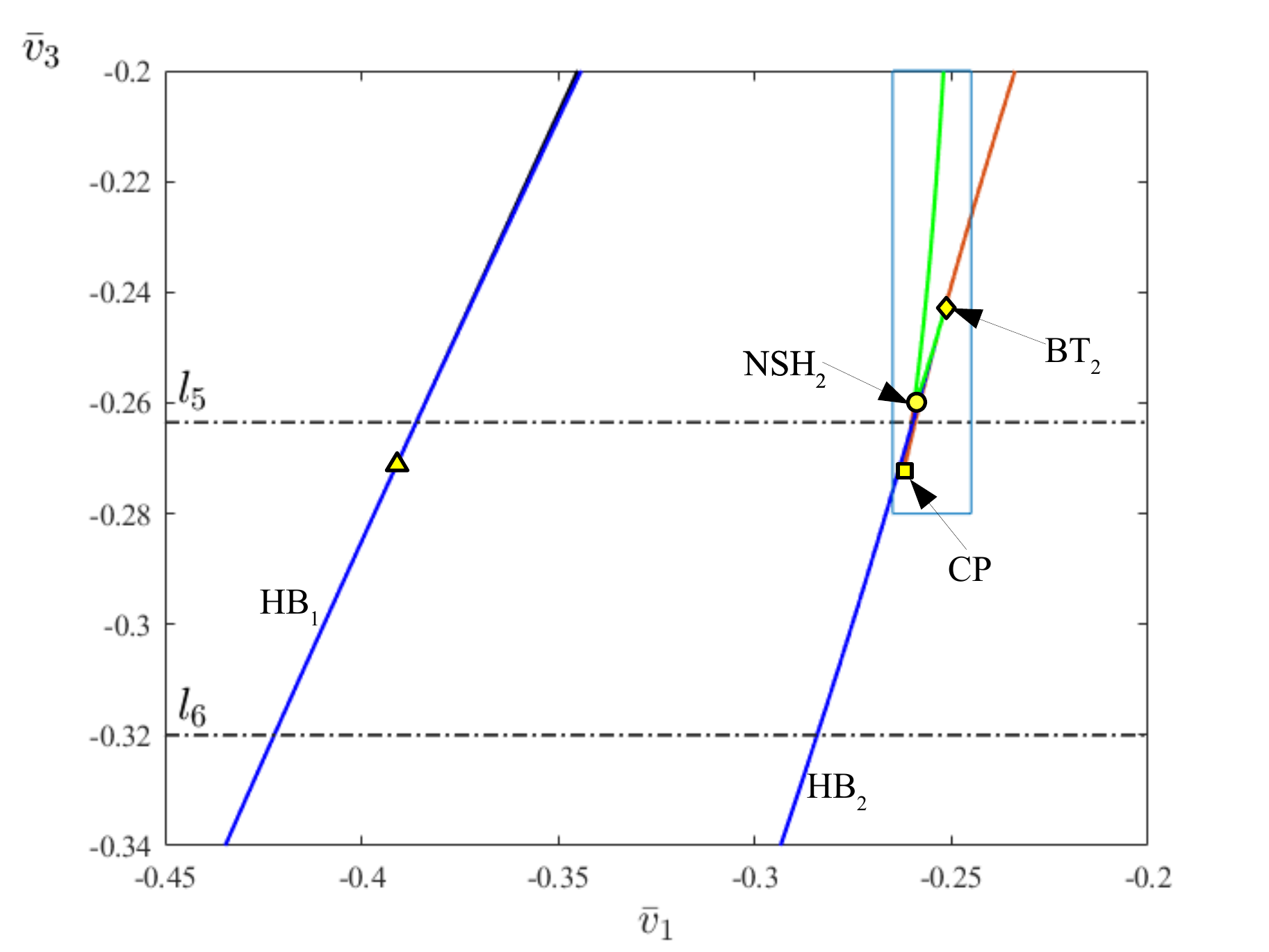}
  \label{fig:11a}
  \end{subfigure}%
  \begin{subfigure}[b]{.56\textwidth}
    \centering
    \caption{}
    \includegraphics[width = \textwidth]{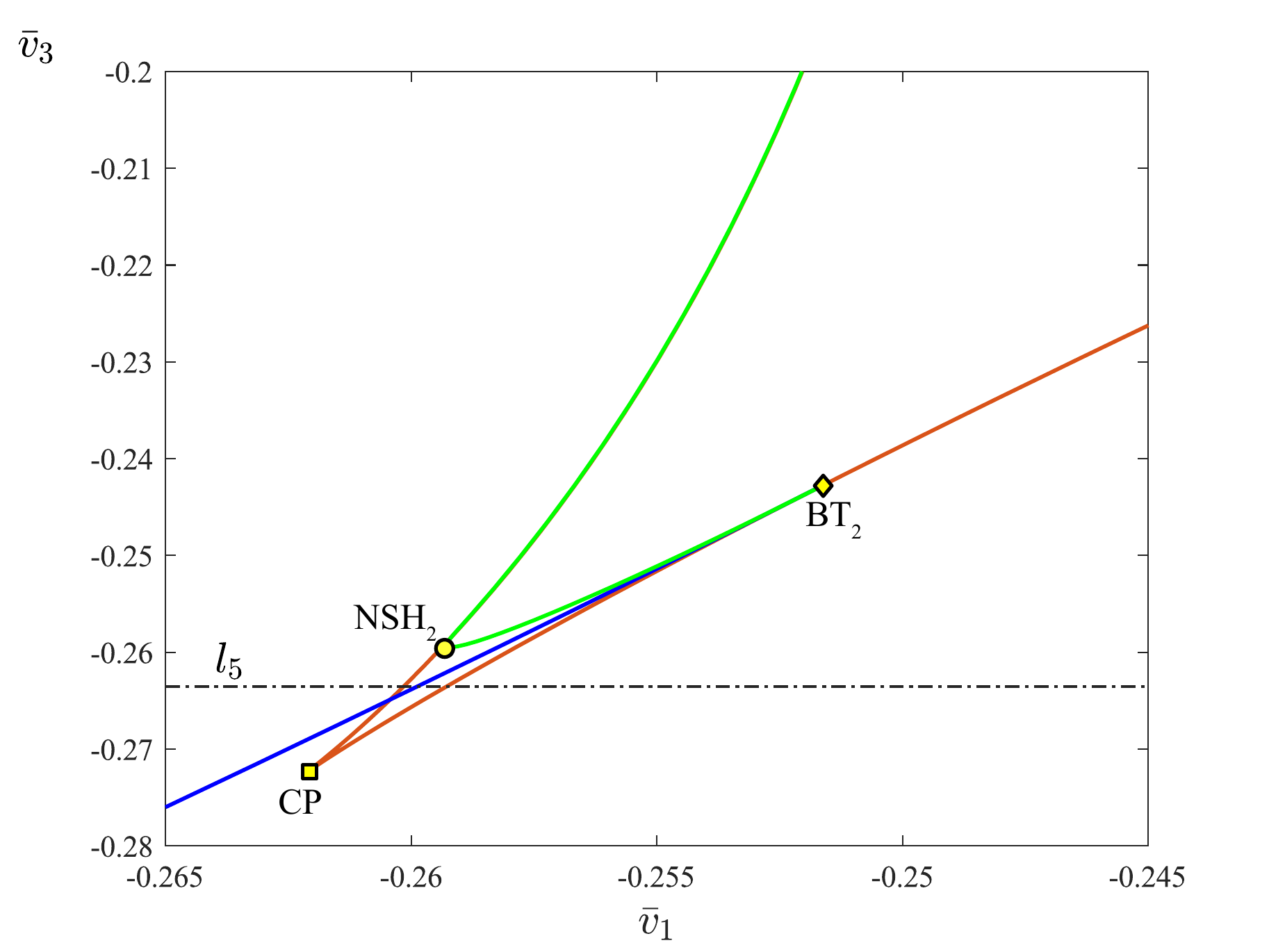}
    \label{fig:11b}
  \end{subfigure}\\%
  \begin{subfigure}[b]{.5\textwidth}
    \centering
    \caption{}
    \includegraphics[width =\textwidth]{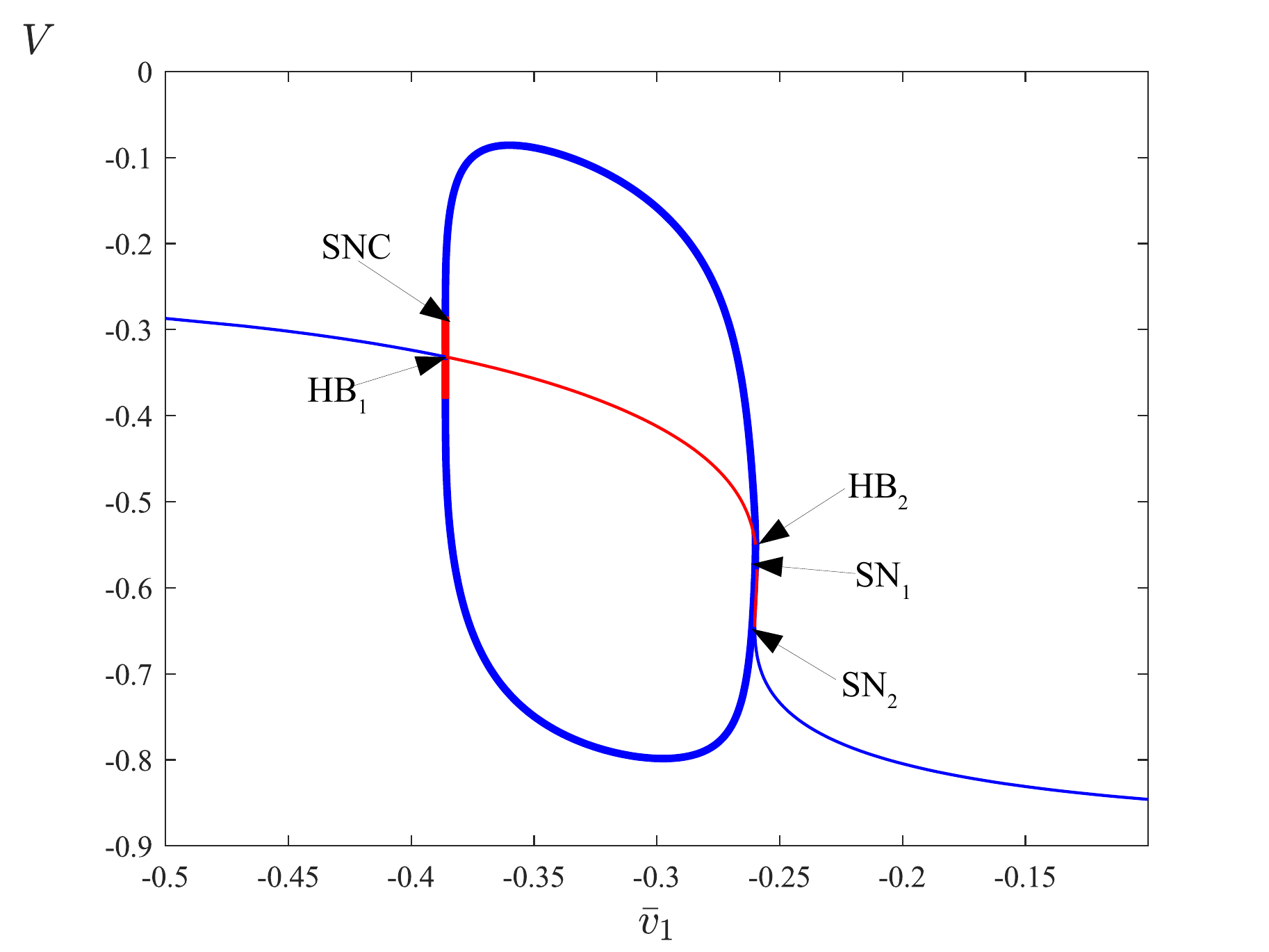}
    \label{fig:11c}
  \end{subfigure}%
  \begin{subfigure}[b]{.5\textwidth}
    \centering
    \caption{}
    \includegraphics[width = \textwidth]{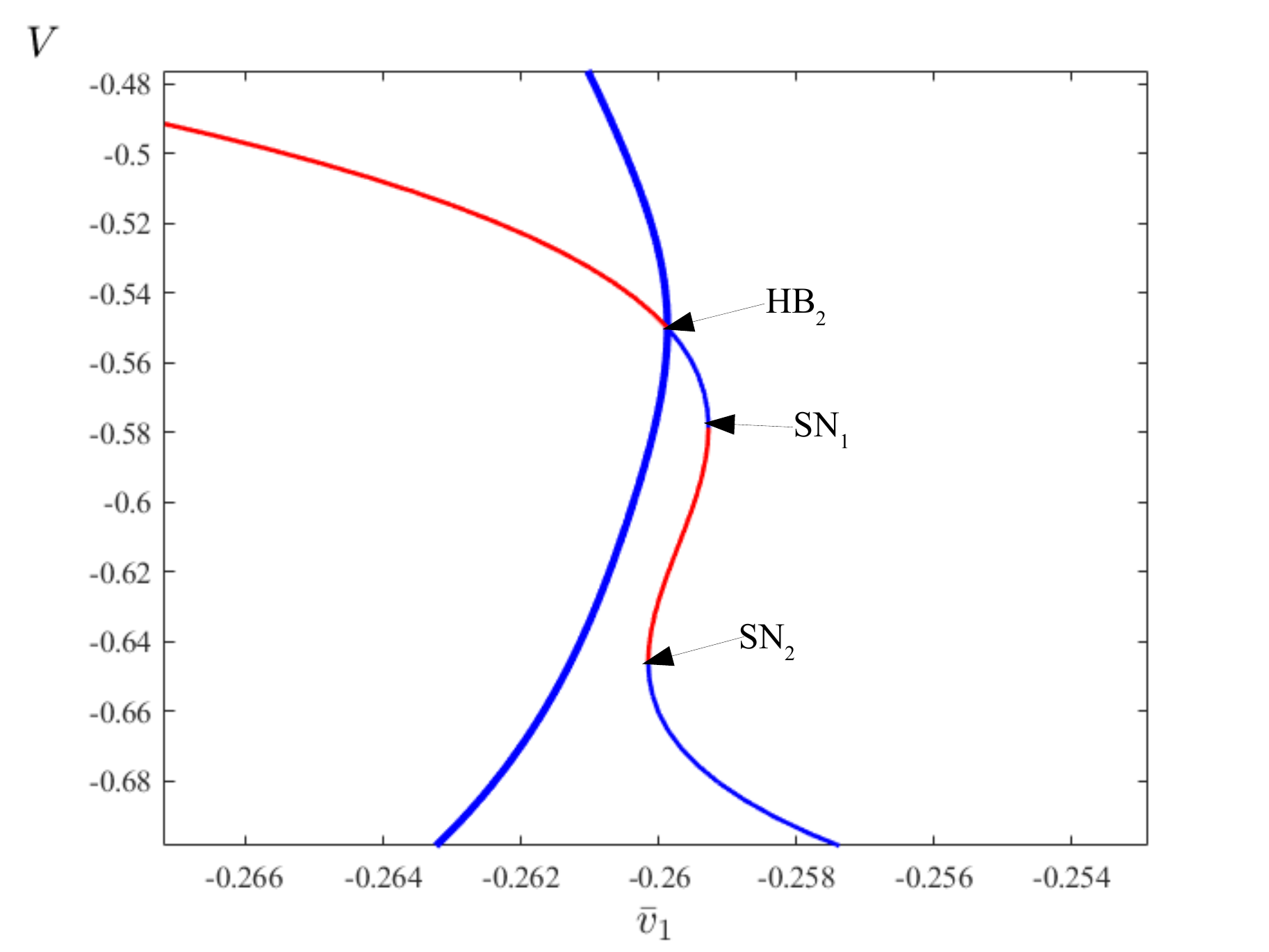}
    \label{fig:11d}
  \end{subfigure}%
 \caption{ (a) An enlargement of Fig.~\ref{fig:7} showing  lines $l_5$ and $l_6$  (b) An enlargement of panel (a). (c)  A one-parameter bifurcation diagram along  $l_5$ with $\bar{v}_3=-0.26$. (d) An enlargement of panel (c). HB: Hopf bifurcation, SN: saddle-node bifurcation, SNC: saddle-node bifurcation of a periodic orbit}
\label{fig:11}
\end{figure}

Upon further decrease to the value of $\bar{v}_3$ a second
Bogdanov-Takens bifurcation, denoted $\text{BT}_2$,
occurs on the $\text{SN}_1$ locus at $\bar{v}_3 \approx -0.2429$ (see Fig.~\ref{fig:11b}). This generates loci of homoclinic and supercritical Hopf bifurcations. The homoclinic locus terminates nearby at another $\text{NSH}_2$ bifurcation where the SNIC locus reverts to a locus of saddle-node bifurcations.
The slice $l_5$, Fig.~\ref{fig:11c}, is below these two codimension-two points. Here the system exhibits Type II excitability as stable oscillations are created at the Hopf bifurcation. This shows that the transition between Type I and Type II excitability for the parameter regime we have considered is governed by the Bogdanov-Takens bifurcation $\text{BT}_2$, and this is in agreement with the result in \citep{Zhao2017TransitionsAutapse} where the authors studied bifurcation mechanisms induced by autapse in the Morris-Lecar model. 

\begin{figure}
      \centering
      \includegraphics[width =0.6\textwidth]{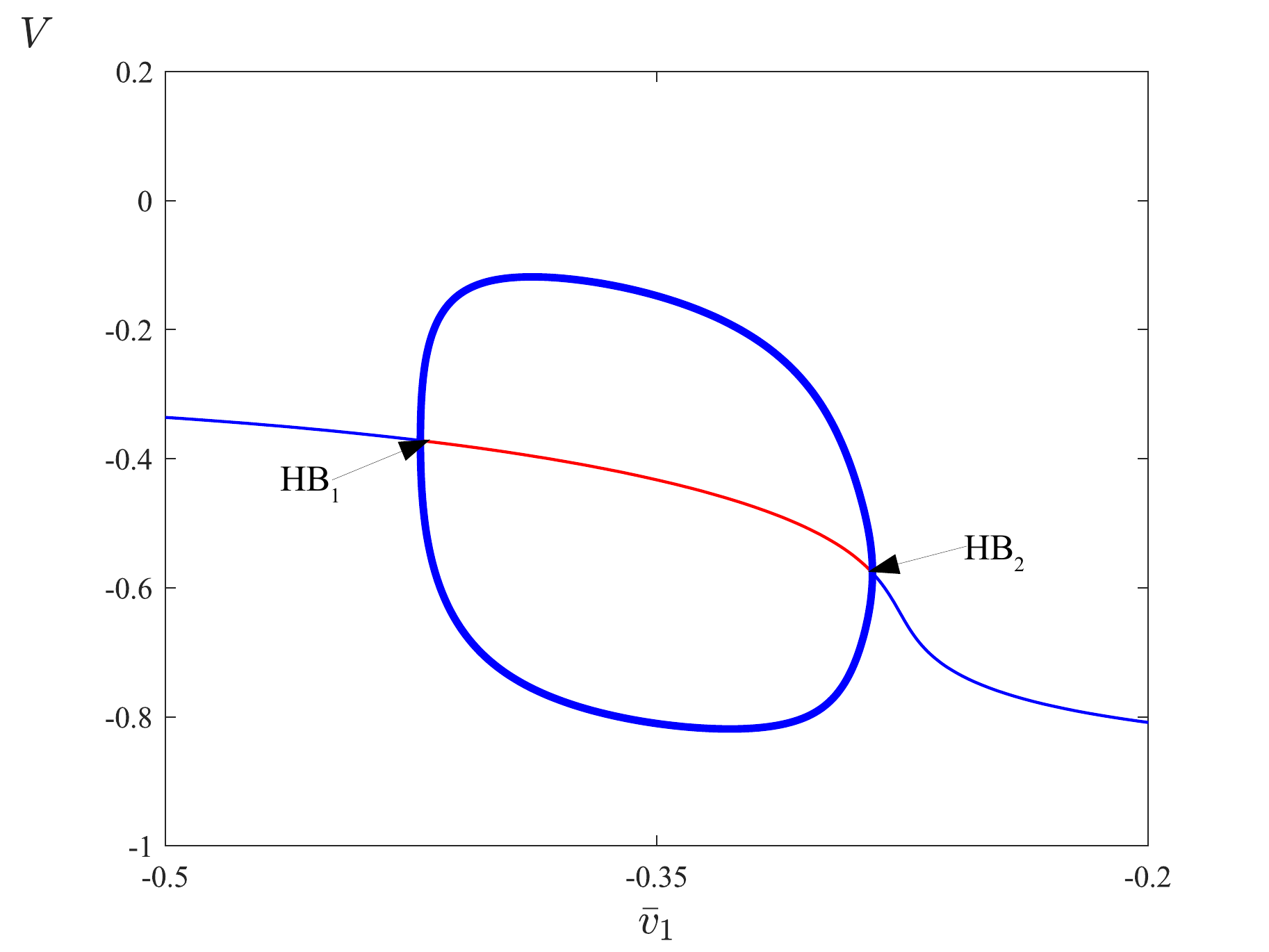}
      \caption{ A one-parameter bifurcation diagram along $l_6$ with $\bar{v}_3=-0.32$ (shown in Figs.~\ref{fig:7} and \ref{fig:11a}). HB: Hopf bifurcation}
      \label{fig:12}
  \end{figure}

Finally, as $\bar{v}_3$ is decreased further the Hopf locus
$\text{HB}_1$ changes from subcritical to supercritical at a generalised Hopf bifurcation at $\bar{v}_3 \approx -0.2708$ and the saddle-node loci $\text{SN}_1$ and $\text{SN}_2$ collide and annihilate in a cusp bifurcation $\text{CP}$ at $\bar{v}_3 \approx -0.2727$. Below these two codimension-two points the only bifurcations that remain
are two supercritical Hopf bifurcations. The slice $l_6$, Fig.~\ref{fig:12}, shows a typical bifurcation diagram. Here the excitability is Type II and there is no bistability.

\newpage
\section{Conclusion}\label{sec:conclusion}
In this paper we have studied a pacemaker model of SMCs where the interactions between ion fluxes, in particular $\text{Ca}^{2+}$ and $\text{K}^{+}$, results in spontaneous oscillations. We established that both $\text{Ca}^{2+}$ and $\text{K}^{+}$ currents are required for the pacemaker activity. Upon varying the voltage associated with the opening of half the  $\text{K}^{+}$ channels, $v_1$, the full three-dimensional model exhibits various dynamical features observed in the conventional models for excitable cells. With the aid of bifurcation diagrams, we showed we show that the reduced two-dimensional model preserves the dynamical properties of the full model qualitatively.

The main motivation of this work was to understand the types of excitability exhibited by the pacemaker model. We showed that the model can be of Type I or Type II excitability depending on how parameters are varied. In particular we determined the bifurcation structure of the $(\bar{v}_1,\bar{v}_3)$-parameter plane to show transitions between the two types of excitability. We found that, as in \cite{Tsumoto2006BifurcationsModel} which used different parameters including non-zero external current, a Bogdanov-Takens bifurcation demarcates the transition between Type I and Type II excitability.

We also revealed that the biologically important parameter $\bar{v}_1$ affects the type of excitability and nature of the oscillations more generally. The results of the model agree with experimental observations on pacemaker behaviour of smooth muscle cells \citep{Meyer1983PacemakerMuscle, Meyer1988IsCells, Harder1984Pressure-DependentArtery, Segal1989ConductionCoupling}  and  neural cells \citep{Connor1985NeuralRhythmicity,Ramirez2004PacemakerView}. 

It is hoped the results may find application in models and experimental studies of physiological and pathophysiological responses in muscle cells. Certainly the observation that the dynamics of SMCs are particularly sensitive to parameter values has been utilised pharmacologically in therapeutics \citep{Droogmans1989ElectromechanicalMuscle,Pogatsa1994AlteredDiabetes}.

Our analysis concerned a single SMC, however SMCs are interconnected through gap junctions and action potentials can propagate between them. It remains to analyse the spatiotemporal behaviour of coupled pacemaker SMCs. Some experimental and computational studies of SMCs have shown that voltage-dependent inward $\text{Na}^{+}$ current is important in EMC activity \citep{Berra-Romani2005TTX-sensitiveCells,Ulyanova2018Voltage-dependentArterioles}, in future work we will incorporate the $\text{Na}^{+}$ current into our model to study its effect on pacemaker dynamics of SMCs.

\acknowledgement

We thank Prof. Hinke M. Osinga (University of Auckland, New Zealand) for the support provided and useful discussion. 

\newpage

\bibliography{sample}

\end{document}